\documentclass[10pt]{amsart}
\usepackage{verbatim}
\usepackage{eucal,url,amssymb,stmaryrd,enumerate,amscd}
\usepackage{amsmath}
\usepackage[pagebackref,colorlinks=true]{hyperref}
\usepackage{amsfonts}
\usepackage{amsthm}

\numberwithin{equation}{section}
\newtheorem{thrm}{Theorem}[section]

\newtheorem{prop}[thrm]{Proposition}

\newtheorem{dfn}[thrm]{Definition}

\newtheorem{rmrk}[thrm]{Remark}

\newtheorem{conv}[thrm]{Convention}
 1  

\begin{document}

\begin{abstract}
We construct explicit left invariant quaternionic contact
structures on Lie groups with zero and non-zero torsion, and  with
non-vanishing quaternionic contact conformal curvature tensor,
thus showing the existence of quaternionic contact manifolds not
locally quaternionic contact conformal to the quaternionic sphere.
We present a left invariant quaternionic contact structure on a
seven dimensional non-nilpotent Lie group, and show that this
structure is locally quaternionic contact conformal to the flat
quaternionic contact structure on the quaternionic Heisenberg
group. On the product of a seven dimensional Lie group, equipped
with a quaternionic contact structure, with the real line we
determine explicit complete quaternionic K\"ahler metrics and
$Spin(7)$-holonomy metrics which seem to be new. We give explicit
complete non-compact eight dimensional { almost quaternion
hermitian} manifolds with closed fundamental four form which are
not quaternionic K\"ahler.
\end{abstract}

\keywords{quaternionic contact structures, Einstein structures, qc conformal flatness, qc conformal curvature,
quaternionic K\"ahler structures, Spin(7)-holonomy metrics, quaternionic K\"ahler and hyper K\"ahler metrics}
\subjclass{58G30, 53C17}

\title[Quaternionic Contact Structures and Special Holonomy Metrics]
{Explicit Quaternionic Contact Structures and
Metrics with Special Holonomy }
\date{\today }
\thanks{This
work has been partially funded by grant MCINN (Spain) MTM2008-06540-C02-01/02}
\author{L.C. de Andr\'es}
\address[Luis C. de Andr\'es, Marisa Fern\'andez, Jos\'e A. Santisteban]{Universidad del Pa%
\'{\i}s Vasco\\
Facultad de Ciencia y Tecnolog\'{\i}a, Departamento de Matem\'aticas\\
Apartado 644, 48080 Bilbao\\
Spain} \email{luisc.deandres@ehu.es} \email{marisa.fernandez@ehu.es} \email{joseba.santisteban@ehu.es}
\author{M. Fern\'andez}
\author{S. Ivanov}
\address[Stefan Ivanov]{University of Sofia, Faculty of Mathematics and
Informatics, blvd. James Bourchier 5, 1164, Sofia, Bulgaria}
\address{and Department of Mathematics,
University of Pennsylvania, Philadelphia, PA 19104-6395}
\email{ivanovsp@fmi.uni-sofia.bg}
\author{J.A. Santisteban}
\author{L. Ugarte}
\address[Luis Ugarte]{Departamento de Matem\'aticas\,-\,I.U.M.A.\\
Universidad de Zaragoza\\
Campus Plaza San Francisco\\
50009 Zaragoza, Spain} \email{ugarte@unizar.es}
\author{D. Vassilev}
\address[Dimiter Vassilev]{ Department of Mathematics and Statistics\\
University of New Mexico\\
Albuquerque, New Mexico, 87131-0001} \email{vassilev@math.unm.edu} \maketitle \tableofcontents

\setcounter{tocdepth}{2}

\section{Introduction}

It is well known that the sphere at infinity of a  non-compact symmetric space $M$ of rank one carries a
natural Carnot-Carath\'eodory structure (see \cite{M,P}).
Quaternionic contact  structures were introduced by Biquard in \cite{Biq1,Biq2}, and they appear naturally
as the conformal boundary at infinity of the quaternionic hyperbolic space. Such structures are also relevant
for the quaternionic contact Yamabe problem which is naturally connected with the extremals and the best
constant in an associated Sobolev-type (Folland-Stein \cite{FS}) embedding on the quaternionic Heisenberg group
\cite{Wei,IMV,IMV1}.

A quaternionic contact structure (\emph{qc structure}) on a real $(4n+3)$-dimensional
manifold $M$, following Biquard, is a codimension three distribution $H$ locally given as the kernel of a
$\mathbb{R}^3$-valued
$1$-form $\eta=(\eta_1,\eta_2,\eta_3)$, such that, the three $2$-forms $%
d\eta_i|_H$ are the fundamental forms of a quaternionic structure on $H$.
The 1-form $\eta$ is determined up to a conformal factor and the action of $SO(3)$ on
$\mathbb{R}^3$, and therefore $H$ is equipped with a conformal class $[g]$ of Riemannian metrics.
The transformations preserving a given quaternionic contact structure $\eta$, i.e.
$\bar\eta=\mu\Psi\eta$ for a positive smooth function $\mu$ and a non-constant$ SO(3)$
matrix $\Psi$ are called \emph{quaternionic contact conformal (qc conformal for short) transformations}.
If the function $\mu$ is constant we
have \emph{quaternionic contact homothetic (qc-homothetic) transformations}.
To every metric in the fixed conformal class $[g]$ on $H$ one can associate a linear connection
preserving the qc structure, see \cite{Biq1}, which we shall call the Biquard connection. The Biquard connection is
invariant under qc homothetic transformations but changes in a non-trivial way under qc conformal transformations.

The quaternionic Heisenberg group $\boldsymbol{G\,(\mathbb{H})}$ with its standard left-invariant qc structure
is the unique (up to a $SO(3)$-action) example of a qc structure with flat Biquard connection  \cite{IV}. The
quaternionic Cayley transform is a quaternionic contact conformal equivalence between the standard 3-Sasakian
structure on the $(4n+3)$-dimensional sphere $S^{4n+3}$ minus a point and the flat qc structure on
$\boldsymbol{G\,(\mathbb{H})}$ \cite{IMV}. All qc structures locally qc conformal to
$\boldsymbol{G\,(\mathbb{H})}$ and $S^{4n+3}$ are characterized in \cite{IV} by the vanishing of a tensor
invariant, the qc-conformal curvature $W^{qc}$ defined in terms of the curvature and torsion of the Biquard
connection.

Examples of qc manifolds  arising from quaternionic K\"ahler deformations are given in \cite{Biq1,Biq2,D1}. A
totally umbilic hypersurface of a quaternionic K\"ahler or hyperK\"ahler manifold carries such a structure. A
basic example is provided by any 3-Sasakian manifold which can be defined as a $(4n+3)$-dimensional Riemannian
manifold whose Riemannian cone is a hyperK\"ahler manifold. It was shown in \cite{IMV} that the torsion
endomorphism of the Biquard connection is the obstruction for a given qc-structure to be locally qc homothetic
to a 3-Sasakian structure provided the scalar curvature of the Biquard connection is positive. Duchemin shows
\cite{D1} that for any qc manifold there exists a quaternionic K\"ahler manifold such that the
qc manifold is realized as a hypersurface. However, the embedding in his construction is not isometric and it
is difficult to write an explicit expression of the quaternionic K\"ahler metric except the 3-Sasakian case
where the cone metric is hyperK\"ahler.

One purpose of this paper is to find new explicit examples of   qc structures. We
construct explicit left invariant qc structures on seven dimensional Lie groups with zero and non-zero torsion
of the Biquard connection for which the qc-conformal curvature tensor does not vanish, $W^{qc}\not=0$ thus
showing the existence of qc manifolds not locally qc conformal to the quaternionic Heisenberg group
$\boldsymbol{G\,(\mathbb{H})}$. We present a left invariant qc strucutre with zero torsion of the Biquard
connection on a seven dimensional non-nilpotent Lie group $G_1$. Surprisingly, we obtain that this qc
structure is locally qc conformal to the flat qc structure on the two-step nilpotent quaternionic Heisenberg group
$\boldsymbol{G\,(\mathbb{H})}$ showing that the qc conformal curvature
is zero and applying the main result in \cite{IV}. Consequently, this fact yields the existence of a local
function $\mu$ such that the qc conformal transformation $\bar\eta=\mu\eta$ preserves  the vanishing of the
torsion of the Biquard connection.

The second goal of the paper is to construct explicit quaternionic K\"ahler and $Spin(7)$-holonomy metrics, i.e
metrics with holonomy  $Sp(n)Sp(1)$ and $Spin(7)$, respectively, on a product of a qc manifold with a real
line. We generalize the notion of a qc structure, namely, we define $Sp(n)Sp(1)$-hypo structures on a
$(4n+3)$-dimensional manifold as structures induced on a hypersurface of a quaternionic K\"ahler manifold.
We present explicit complete non-compact quaternionic K\"ahler metrics and
$Spin(7)$-holonomy metrics on the product of the locally qc conformally flat quaternionic contact structure on
the seven dimensional Lie group $G_1$ with the real line some of which seem to be new, see Section 5 and 6.

It is well known that in dimension eight { an almost quaternion
hermitian} structure with closed fundamental  four form is not
necessarily quaternionic K\"ahler \cite{Sw}. This fact was
confirmed by Salamon constructing in \cite{Sal} a compact example
of { an almost quaternion hermitian} manifold with closed
fundamental four form which is not Einstein, and therefore it is
not a quaternionic K\"ahler. We give a three parameter family of
explicit complete non-compact eight dimensional { almost
quaternion hermitian manifolds} with closed fundamental four form
which are not quaternionic K\"ahler. We also check that these
examples are not Einstein as well.

To the best of our knowledge there is not { known example of  an
almost quaternion hermitian}  eight dimensional manifold with
closed fundamental four form which is Einstein but not
quaternionic K\"ahler.

In dimension four, we recover some of the known hyper K\"ahler metrics known as gravitational instantons
(Bianchi-type metrics). Furthermore, we give explicit hyper-symplectic  (hyper para K\"ahler) metrics of
signature (2,2).  A hyper symplectic structure in dimension four underlines an anti-self-dual Ricci-flat
neutral metric. For this reason such structures  have been used in string theory \cite{OV,hul,JR,Bar,Hull,CHO}
and integrable systems \cite{D12,BM,DW}. Our construction gives a kind of duality between hyper K\"ahler and
hyper para K\"ahler structures in dimension four.

\begin{conv}
\label{conven} \hfill\break\vspace{-15pt}
\begin{enumerate}
\item[a)] We shall use $X,Y,Z,U$ to denote horizontal vector fields, i.e. $%
X,Y,Z,U\in H$;
\item[b)] $\{e_1,\dots,e_{4n}\}$ denotes a  local orthonormal basis of the horizontal
space $H$;
\item[c)] The summation convention over repeated vectors from the basis $%
\{e_1,\dots,e_{4n}\}$ is used. For example, the formula $k=P(e_b,e_a,e_a,e_b) $ means
$k=\sum_{a,b=1}^{4n}P(e_b,e_a,e_a,e_b). $
\item[d)] The triple $(i,j,k)$ denotes any cyclic permutation of $(1,2,3)$.
\item[e)] $s$  will be any number from the set $\{1,2,3\}, \quad
s\in\{1,2,3\}$.
\end{enumerate}
\end{conv}

\textbf{Acknowledgments} We thank  Charles Boyer for useful conversions leading to Remark \ref{r:CBoyer} and
Dieter Lorenz-Petzold for several comments on the Bianchi type solutions.

The research was initiated during the visit of the third author to the Abdus Salam ICTP, Trieste as a Senior
Associate, Fall 2008. He also thanks ICTP for providing the support and an excellent research environment. S.I.
is partially supported by the Contract 082/2009 with the University of Sofia `St.Kl.Ohridski'. S.I and D.V. are
partially supported by Contract ``Idei", DO 02-257/18.12.2008. This work has been also partially supported
through grant MCINN (Spain) MTM2008-06540-C02-01/02.

\section{Quaternionic contact manifolds}

In this section we will briefly review the basic notions of quaternionic contact geometry and recall some
results from \cite{Biq1}, \cite{IMV} and \cite{IV} which we will use in this paper.

\subsection{qc structures and the Biquard connection}

A quaternionic contact (qc) manifold $(M, g, \mathbb{Q})$ is a $4n+3$-dimensional manifold $M$ with a
codimension three distribution $H$ satisfying
\begin{enumerate}
\item[i)] { $H$ has an $Sp(n)Sp(1)$ structure, that is it is
equiped with a Riemannian metric $g$ and a rank-three bundle
$\mathbb Q$ consisting of (1,1)-tensors on $H$ locally generated
by three almost complex structures $I_1,I_2,I_3$ on $H$ satisfying
the identities of the imaginary unit quaternions,
$I_1I_2=-I_2I_1=I_3, \quad I_1I_2I_3=-id_{|_H}$ which are
hermitian compatible with the metric $g(I_s.,I_s.)=g(.,.), s =
1,2,3$, i.e., $H$ has an almost quaternion hermitian structure.}



\item[ii)]  $H$ is locally given as the kernel of a 1-form $\eta=(\eta_1,\eta_2,\eta_3)$ with
values in $\mathbb{R}^3$ and the following compatibility condition holds $%
\qquad 2g(I_sX,Y)\ =\ d\eta_s(X,Y), \quad s=1,2,3, \quad X,Y\in H.$
\end{enumerate}

A special phenomena here, noted in \cite{Biq1}, is that the
contact form $\eta$ determines the  quaternionic structure and the
metric on the horizontal distribution in a unique way.

Correspondingly, given a quaternionic contact manifold we shall denote with $%
\eta$ any associated contact form. The associated contact form is
determined up to  an $SO(3)$-action, namely if $\Psi\in SO(3)$
then $\Psi\eta$ is again a contact form satisfying the above
compatibility condition (rotating also the almost complex
structures). On the other hand, if we consider the conformal class
$[g]$ on $H$, the associated contact forms are determined up to a
multiplication with a positive conformal factor $\mu$ and an
$SO(3)$-action, namely if $\Psi\in SO(3)$ then $\mu\Psi\eta$ is a
contact form associated
with a metric in the conformal class $[g]$ on $H$. A qc manifold $(M, \bar g,\mathbb{Q} )$
is called qc conformal to $(M, g,%
\mathbb{Q} )$ if $\bar g\in [g]$. In that case, if $\bar\eta$ is a
corresponding associated $1$-form with complex structures $\bar I_s$, $%
s=1,2,3,$ we have $\bar\eta\ =\ \mu\, \Psi\,\eta$ for some $\Psi\in SO(3)$ with smooth functions as entries and
a positive function $\mu$. In particular, starting with a qc manifold $(M, \eta)$ and defining $\bar\eta\ =\
\mu\, \eta$ we obtain a qc manifold $(M, \bar\eta)$ qc conformal to the original one.

If the first { Pontryagin} class of $M$ vanishes then the 2-sphere bundle of $%
\mathbb{R}^3$-valued 1-forms is trivial \cite{AK}, i.e. there is a globally
defined form $\eta$ that
anihilates $H$, we denote the corresponding qc manifold $(M,\eta)$. In this case the 2-sphere of associated
almost complex structures is also globally defined on $H$.

Any endomorphism $\Psi$ of $H$ decomposes with respect to the quaternionic structure $(\mathbb{Q},g)$ uniquely
into $Sp(n)$-invariant parts as follows
\hspace{2mm} 
$\Psi=\Psi^{+++}+\Psi^{+--}+\Psi^{-+-}+\Psi^{--+}, $ 
\hspace{2mm} where $\Psi^{+++}$ commutes with all three $I_i$, $\Psi^{+--}$
commutes with $I_1$ and anti-commutes with the other two and etc. The two $%
Sp(n)Sp(1)$-invariant components are given by $\Psi_{[3]}=\Psi^{+++}, \quad
\Psi_{[-1]}=\Psi^{+--}+\Psi^{-+-}+\Psi^{--+}. $ Denoting the corresponding
(0,2) tensor via $g$ by the same letter one sees that the $Sp(n)Sp(1)$%
-invariant components are the projections on the eigenspaces of the Casimir operator $\dag = I_1\otimes I_1 +
I_2\otimes I_2 + I_3\otimes I_3 $
corresponding, respectively, to the eigenvalues $3$ and $-1$, see \cite{CSal}%
. If $n=1$ then the space of symmetric endomorphisms commuting with all $%
I_i, i=1,2,3$ is 1-dimensional, i.e. the [3]-component of any symmetric
endomorphism $\Psi$ on $H$ is proportional to the identity, $\Psi_{[3]}=%
\frac{tr\, (\Psi)}{4}Id_{|H}$.

On a quaternionic contact manifold there exists a canonical connection defined in \cite{Biq1} when the
dimension $(4n+3)>7$, and in \cite{D} in the 7-dimensional case.

\begin{thrm}
\cite{Biq1}\label{biqcon} {Let $(M, g,\mathbb{Q})$ be a quaternionic contact manifold} of dimension $4n+3>7$
and a fixed metric $g$ on $H$ in the conformal class $[g]$. Then there exists a unique connection $\nabla$ with
torsion $T$ on $M^{4n+3}$ and a unique supplementary subspace $V$ to $H$ in $%
TM$, such that:

\begin{enumerate}
\item[i)] $\nabla$ preserves the decomposition $H\oplus V$ and the
metric $g$; \item[ii)] for $X,Y\in H$, one has
$T(X,Y)=-[X,Y]_{|V}$; \item[iii)] { $\nabla$ preserves the
$Sp(n)Sp(1)$ structure on $H$, i.e. $\nabla g=0,  \nabla\sigma
\in\Gamma(\mathbb Q)$ for a section $\sigma\in\Gamma(\mathbb Q)$;}
\item[iv)] for $\xi\in V$, the endomorphism $T(\xi,.)_{|H}$ of $H$ lies in $%
(sp(n)\oplus sp(1))^{\bot}\subset gl(4n)$;
\item[v)] the connection on $V$ is induced by the natural identification $\varphi$ of
$V$ with the subspace $sp(1)$ of the endomorphisms of $H$, i.e. $\nabla\varphi=0$.
\end{enumerate}
\end{thrm}
In iv), the inner product $<,>$ of $End(H)$ is given by $<A,B> = { \sum_{i=1}^{4n} g(A(e_i),B(e_i)),}$
for $A, B \in End(H)$.

We shall call the above connection \emph{the Biquard connection}.
Biquard \cite{Biq1} also described the supplementary subspace $V$, namely, {%
locally }$V$ is generated by vector fields $\{\xi_1,\xi_2,\xi_3\}$, such that
\begin{equation}  \label{bi1}
\begin{aligned} \eta_s(\xi_k)=\delta_{sk}, \qquad (\xi_s\lrcorner
d\eta_s)_{|H}=0,\\ (\xi_s\lrcorner d\eta_k)_{|H}=-(\xi_k\lrcorner d\eta_s)_{|H}, \end{aligned}
\end{equation}
where $\lrcorner$ denotes the interior multiplication. The vector fields $\xi_1,\xi_2,\xi_3$ are called Reeb
vector fields.

{ If the dimension of $M$ is seven, { there might be no vector
fields satisfying \eqref{bi1}}.  Duchemin shows in \cite{D} that
if we assume, in addition, the existence of Reeb vector fields as
in \eqref{bi1}, then Theorem~\ref{biqcon}
holds. Henceforth, by a qc structure in dimension $7$ we shall 
mean a qc structure satisfying \eqref{bi1}.}

Notice that equations \eqref{bi1} are invariant under the natural $SO(3)$
action. Using the triple of Reeb vector fields we extend $g$ to a metric on $%
M$ by requiring 
$span\{\xi_1,\xi_2,\xi_3\}=V\perp H \text{ and } g(\xi_s,\xi_k)=\delta_{sk}.
$ 
\hspace{2mm} \noindent The extended metric does not depend on the action of $%
SO(3)$ on $V$, but it changes in an obvious manner if $\eta$ is multiplied by a conformal factor. Clearly, the
Biquard connection preserves the extended metric on $TM, \nabla g=0$.

The covariant derivative of the qc structure with respect to the Biquard connection and the covariant
derivative of the distribution $V $ are given by
\begin{equation}  \label{der}
\nabla I_i=-\alpha_j\otimes I_k+\alpha_k\otimes I_j,\qquad
\nabla\xi_i=-\alpha_j\otimes\xi_k+\alpha_k\otimes\xi_j,
\end{equation}
where the $sp(1)$-connection 1-forms $\alpha_s$ on $H$ are given by \cite%
{Biq1}
\begin{gather}  \label{coneforms}
\alpha_i(X)=d\eta_k(\xi_j,X)=-d\eta_j(\xi_k,X), \quad X\in H, \quad \xi_i\in V,
\end{gather}
while the $sp(1)$-connection 1-forms $\alpha_s$ on the vertical space $V$ are calculated in \cite{IMV}
\begin{gather}  \label{coneform1}
\alpha_i(\xi_s)\ =\ d\eta_s(\xi_j,\xi_k) -\ \delta_{is}\left(\frac{S}2\ +\ \frac12\,\left(\,
d\eta_1(\xi_2,\xi_3)\ +\ d\eta_2(\xi_3,\xi_1)\ + \ d\eta_3(\xi_1,\xi_2)\right)\right),
\end{gather}
where $s\in\{1,2,3\}$ and $S$ is the \emph{normalized} qc scalar curvature defined
below in \eqref{qscs}. The vanishing of the $sp(1)$-connection 1-forms on $%
H $ implies the vanishing of the torsion endomorphism of the Biquard connection {(see \cite{IMV})}.

The fundamental 2-forms $\omega_i, i=1,2,3$ \cite{Biq1}  are
defined by
$2\omega_{i|H}\ =\ \, d\eta_{i|H},\quad \xi\lrcorner\omega_i=0,\quad \xi\in V.
$
The properties of the Biquard connection are encoded in the properties of
the torsion endomorphism $T_{\xi}=T(\xi,\cdot) : H\rightarrow H, \quad \xi\in V$%
. Decomposing the endomorphism $T_{\xi}\in(sp(n)+sp(1))^{\perp}$ into its symmetric part $T^0_{\xi}$ and
skew-symmetric part $b_{\xi}, T_{\xi}=T^0_{\xi} + b_{\xi} $, O. Biquard in \cite{Biq1} shows that the torsion
$T_{\xi}$ is completely trace-free, $tr\, T_{\xi}=tr\, T_{\xi}\circ
I=0, \quad I\in {\mathbb{Q}}$, its symmetric part has the properties $%
T^0_{\xi_i}I_i=-I_iT^0_{\xi_i}\quad I_2(T^0_{\xi_2})^{+--}=I_1(T^0_{\xi_1})^{-+-},\quad
I_3(T^0_{\xi_3})^{-+-}=I_2(T^0_{\xi_2})^{--+},\quad I_1(T^0_{\xi_1})^{--+}=I_3(T^0_{\xi_3})^{+--} $ and the
skew-symmetric part can be represented as $b_{\xi_i}=I_iu$, where $u$ is a traceless symmetric (1,1)-tensor on
$H$ which commutes with $I_1,I_2,I_3$. If $n=1$ then the tensor $u$ vanishes identically, $u=0$ and the torsion
is a symmetric tensor, $T_{\xi}=T^0_{\xi}$.

Any 3-Sasakian manifold has zero torsion endomorphism, and
 the converse is true if in addition the normalized qc scalar curvature $S$ (see
\eqref{qscs}) is a positive constant \cite{IMV}. We remind that a
$(4n+3)$-dimensional  Riemannian manifold $(M,g)$ is
called 3-Sasakian if the cone metric $g_c=t^2g+dt^2$ on $C=M\times \mathbb{R}%
^+$ is a hyper K\"ahler metric, namely, it has holonomy contained in $Sp(n+1)$ \cite{BGN}. A 3-Sasakian
manifold of dimension $(4n+3)$ is Einstein with positive Riemannian scalar curvature $(4n+2)(4n+3)$ \cite{Kas}
and if complete it is compact with a finite fundamental group, due to Mayer's theorem (see \cite{BG} for a nice
overview of 3-Sasakian spaces).

\subsection{Torsion and curvature}

Let $R=[\nabla,\nabla]-\nabla_{[\ ,\ ]}$ be the curvature tensor of $\nabla$ and the dimension is $4n+3$. We
denote the curvature tensor of type (0,4) by the same letter, $R(A,B,C,D):=g(R(A,B)C,D)$, $A,B,C,D \in
\Gamma(TM)$. The Ricci $2$-forms and the normalized scalar curvature of the Biquard connection, called
\emph{qc-Ricci forms} and \emph{normalized qc-scalar curvature}, respectively, are defined by
\begin{equation}  \label{qscs}
4n\rho_s(X,Y)=R(X,Y,e_a,I_se_a), \quad 8n(n+2)S=R(e_b,e_a,e_a,e_b).
\end{equation}
The $sp(1)$-part of $R$ is determined by the Ricci 2-forms and the connection 1-forms by
\begin{equation}  \label{sp1curv}
R(A,B,\xi_i,\xi_j)=2\rho_k(A,B)=(d\alpha_k+\alpha_i\wedge\alpha_j)(A,B), \qquad A,B \in \Gamma(TM).
\end{equation}
The structure equations of a qc structure, discovered in \cite{IV1}, read
\begin{gather}  \label{streq}
2\omega_i=d\eta_i+\eta_j\wedge\alpha_k-\eta_k\wedge\alpha_j + S \eta_j\wedge\eta_k
\end{gather}
and the qc structure is 3-Sasakian exactly when
\begin{gather}  \label{streq3}
2\omega_i=d\eta_i-2 \eta_j\wedge\eta_k,
\end{gather}

\noindent for any  cyclic permutation $(i,j,k)$ of $(1,2,3)$. The two $Sp(n)Sp(1)$-invariant trace-free
symmetric 2-tensors $T^0, U$ on $H$
are introduced in \cite{IMV} as follows $T^0(X,Y)\overset{def}{=}%
g((T_{\xi_1}^{0}I_1+T_{\xi_2}^{0}I_2+T_{ \xi_3}^{0}I_3)X,Y)$, $U(X,Y)%
\overset{def}{=}g(uX,Y)$. The tensor $T^0$ belongs to the [-1]-eigenspace while $U$ is in the [3]-eigenspace of
the operator $\dag$, i.e., they have the properties:
\begin{equation}  \label{propt}
\begin{aligned} T^0(X,Y)+T^0(I_1X,I_1Y)+T^0(I_2X,I_2Y)+T^0(I_3X,I_3Y)=0, \\
U(X,Y)=U(I_1X,I_1Y)=U(I_2X,I_2Y)=U(I_3X,I_3Y). \end{aligned}
\end{equation}
In dimension seven $(n=1)$, the tensor $U$ vanishes identically, $U=0$.

We shall need the following identity taken from \cite[Proposition~2.3]{IV}
\begin{equation}  \label{need}
4g(T^0(\xi_s,I_sX),Y)=T^0(X,Y)-T^0(I_sX,I_sY)
\end{equation}

The horizontal Ricci 2-forms can be expressed in terms of the torsion of the Biquard connection \cite{IMV} (see
also \cite{IMV1,IV}). We collect the necessary facts from \cite[Theorem~1.3, Theorem~3.12, Corollary~3.14,
Proposition~4.3 and Proposition~4.4]{IMV} with slight modification presented in \cite{IV}

\begin{thrm}
\cite{IMV}\label{sixtyseven} 
On a $(4n+3)$-dimensional qc manifold $(M,\eta,\mathbb{Q})$ the next formulas hold
\begin{equation}  \label{sixtyfour}
\begin{aligned} 2\rho_s(X,I_sY) \ & =\
-T^0(X,Y)-T^0(I_sX,I_sY)-4U(X,Y)-2Sg(X,Y),\\ T(\xi_{i},\xi_{j})& =-S\xi_{k}-[\xi_{i},\xi_{j}]_{H}.
\end{aligned}
\end{equation}
The vanishing of the trace-free part of the Ricci 2-forms is equivalent to the vanishing of the torsion
endomorphism of the Biquard connection. In this case the vertical distribution is integrable, the
(normalized) qc scalar
curvature $S$ is constant and if $S>0$ then there locally exists an $SO(3)$%
-matrix $\Psi$ with smooth entries depending on an auxiliary parameter such the (local) qc structure
$(\frac{S}2\Psi\eta,\mathbb{Q})$ is 3-Sasakian.
\end{thrm}

If dimension is bigger than seven it turns out that the vanishing
of the torsion endomorphism of the Biquard connection is
equivalent the $4$-form $
\Omega=\omega_1\wedge\omega_1+\omega_2\wedge\omega_2+\omega_3\wedge\omega_3$
to be closed \cite{IV1}.

\subsection{The qc conformal curvature}

The qc conformal curvature tensor $W^{qc}$ introduced  in \cite{IV} is the obstruction for a qc structure to be
locally qc
conformal to the flat structure on the quaternionic Heisenberg group $%
\boldsymbol{G\,(\mathbb{H})}$. In terms of the torsion and curvature of the Biquard connection $W^{qc}$ is
defined in \cite{IV} by
\begin{multline}  \label{qccm}
W^{qc}(X,Y,Z,V) =\frac14\Big[R(X,Y,Z,V)+\sum_{s=1}^3R(I_sX,I_sY,Z,V)\Big] \\
+(g\owedge U) (X,Y,Z,V) + \sum_{s=1}^3(\omega_s\owedge I_sU)(X,Y,Z,V)
-\frac12\sum_{s=1}^3\omega_s(Z,V)\Big[T^0(X,I_sY)-T^0(I_sX,Y)\Bigr] \\
+\frac{S}4\Big [ (g\owedge g)(X,Y,Z,V) +\sum_{s=1}^3 (\omega_s \owedge %
\omega_s)(X,Y,Z,V) \Big ],
\end{multline}
where $I_s U\, (X,Y) = -U (X,I_s Y)$ and $\owedge$ is the Kulkarni-Nomizu product of (0,2) tensors, for
example,
\begin{multline*}
(\omega_s\owedge U)(X,Y,Z,V):=\omega_s(X,Z)U(Y,V)+ \omega_s(Y,V)U(X,Z)-\omega_s(Y,Z)U(X,V)-\omega_s(X,V)U(Y,Z).
\end{multline*}
The main result from \cite{IV} can be stated as follows

\begin{thrm}
\label{main1}\cite{IV} A qc structure on a $(4n+3)$-dimensional smooth manifold is locally quaternionic contact
conformal to the standard flat qc structure on the quaternionic Heisenberg group $\boldsymbol{G\,(\mathbb{H})}$
if and only if the qc conformal curvature vanishes, $W^{qc}=0$. In this case, we call the qc structure a qc
conformally flat structure.
\end{thrm}

Denote $L_0=\frac12T^0+U$. For computational purposes we use the fact established in \cite{IV} that $W^{qc}=0$
exactly when the tensor $WR=0$, where
\begin{multline}  \label{qcwdef1}
WR(X,Y,Z,V)= R(X,Y,Z,V)+ (g\owedge L_0))(X,Y,Z,V)+\sum_{s=1}^3(\omega_s%
\owedge I_sL_0)(X,Y,Z,V) \\
-\frac12\sum_{s=1}^3\Bigl[\omega_s(X,Y)\Bigl\{T^0(Z,I_sV)-T^0(I_sZ,V)\Bigr\} %
+ \omega_s(Z,V)\Bigl\{T^0(X,I_sY)-T^0(I_sX,Y)-4U(X,I_sY)\Bigr\}\Bigr] \\
+\frac{S}4\Big[(g\owedge g)(X,Y,Z,V)+\sum_{s=1}^3\Bigl((\omega_s\owedge%
\omega_s)(X,Y,Z,V) +4\omega_s(X,Y)\omega_s(Z,V)\Bigr) \Big].
\end{multline}
We also recall that as a manifold $\boldsymbol{G\,(\mathbb{H})} \ =\mathbb{H}%
^n\times\text {Im}\, \mathbb{H}$, while the group multiplication is given by $( q^{\prime }, \omega^{\prime })\
=\ (q_o, \omega_o)\circ(q, \omega)\ =\
(q_o\ +\ q, \omega\ +\ \omega_o\ + \ 2\ \text {Im}\ q_o\, \bar q)$, where $%
q,\ q_o\in\mathbb{H}^n$ and $\omega, \omega_o\in \text {Im}\, \mathbb{H}$. The standard flat quaternionic
contact structure is defined by the left-invariant quaternionic contact form $\tilde\Theta\ =\
(\tilde\Theta_1,\ \tilde\Theta_2, \ \tilde\Theta_3)\ =\ \frac 12\ (d\omega \ - \ q^{\prime
}\cdot d\bar q^{\prime }\ + \ dq^{\prime }\, \cdot\bar q^{\prime })$, where $%
.$ denotes the quaternion multiplication. As a Lie group it can be characterized by the following structure
equations. Denote by ${e^a, 1 \leq a\leq (4n+3)}$ the
basis of the left invariant 1-forms,
{and by $e^{ij}$ the wedge product $e^i \wedge e^j$.} The $(4n+3)$-dimensional
quaternionic Heisenberg Lie algebra is the 2-step nilpotent Lie
algebra defined by:
\begin{equation}  \label{4n+3heis}
\begin{aligned}
& de^a =0, \qquad 1\leq a \leq 4n, \\
& d\eta_1=de^{4n+1}= 2( e^{12} + e^{34}+\cdots + e^{(4n-3)(4n-2)}+ e^{(4n-1)4n} ),
\\
& d\eta_2=de^{4n+2}= 2 (e^{13} +e^{42}+\cdots + e^{(4n-3)(4n-1)}+ e^{4n(4n-2)}),
\\
& d\eta_3 = de^{4n+3}= 2 (e^{14} + e^{23}+\cdots + e^{(4n-3)4n}+ e^{(4n-2)(4n-1)}).
\end{aligned}
\end{equation}


\section{Examples}

In this section we give explicit examples of qc structures in dimension seven satisfying the compatibility
conditions \eqref{bi1}. The first example has zero torsion and is locally qc conformal to the quaternionic
Heisenberg group. The second example has zero torsion while the third is with non-vanishing torsion, and both
are not locally qc conformal to the
quaternionic Heisenberg group. 

Clearly, a qc conformally flat structure is locally qc conformal to a 3-Sasaki structure due to the local qc
conformal equivalence of the standard 3-Sasakian structure on the $4n+3$-dimensional sphere and the
quaternionic Heisenberg group.

\begin{rmrk}\label{r:CBoyer}
We note explicitly that the vanishing of the torsion endomorphism
implies that, locally,  the structure is homothetic to a
3-Sasakian structure if the qc scalar curvature is positive. In
the seven dimensional examples below the qc scalar curvature is a
negative constant. In that respect, as pointed by Charles Boyer,
there are no compact invariant with respect to translations
3-Sasakian Lie groups of dimension seven.
\end{rmrk}

\subsection{Zero torsion qc-flat-Example 1}

Denote $\{\tilde e^l,\ 1 \leq l\leq 7\}$ the basis of the left invariant 1-forms and consider the simply
connected Lie group with Lie algebra $\widetilde{L}_1$ defined by the following equations:
\begin{equation}\label{exb1}
\begin{aligned}
&d\tilde e^1=0,\quad d\tilde e^2=\tilde e^{34}, \quad d\tilde e^3=-\tilde e^{24},\quad d\tilde e^4=\tilde
e^{23},\quad
d\tilde e^5=-2\tilde e^{14}-2\tilde e^{23}+\tilde e^{15}+\tilde e^{26}-\tilde e^{37},\\
&d\tilde e^6=-2\tilde e^{13}-2\tilde e^{42}+\tilde e^{16}-\tilde e^{25}+\tilde e^{47},\quad d\tilde
e^7=-2\tilde e^{12}-2\tilde e^{34}+\tilde e^{17}+\tilde e^{35}-\tilde e^{46}.
\end{aligned}
\end{equation}
Let $L_1$ be the Lie algebra  isomorphic to  \eqref{exb1} described by
\begin{equation}  \label{ex11}
\begin{aligned}
&de^1=0,\quad de^2=-e^{12}-2e^{34}-\frac12e^{37}+\frac12e^{46},\\
&de^3=-e^{13}+2e^{24}+\frac12e^{27}-\frac12e^{45},\quad
de^4=-e^{14}-2e^{23}-\frac12e^{26}+\frac12e^{35}\\
&de^5=2e^{12}+2e^{34}-\frac12e^{67},\quad de^6=2e^{13}+2e^{42}+\frac12e^{57},\quad
de^7=2e^{14}+2e^{23}-\frac12e^{56}.
\end{aligned}
\end{equation}
and $e_l, 1 \leq l\leq 7$ be the left invariant vector field dual to the 1-forms ${e^i, 1 \leq i\leq 7}$,
respectively. We define a global qc structure on $L_1$ by setting
\begin{equation}  \label{qc1}
\begin{aligned} &\eta_1=e^5, \quad \eta_2=e^6, \quad \eta_3=e^7, \quad
H=span\{e^1,\dots, e^4\},\\ &\omega_1=e^{12}+e^{34}, \quad \omega_2=e^{13}+e^{42}, \quad
\omega_3=e^{14}+e^{23}. \end{aligned}
\end{equation}
It is straightforward to check from \eqref{ex11} that the vector fields $
\xi_1=e_5$, $\xi_2=e_6$, $\xi_3=e_7 $ 
satisfy the Duchemin compatibility conditions \eqref{bi1} and therefore the Biquard connection exists and
$\xi_s$ are the Reeb vector fields.

\begin{thrm}
\label{m1} Let $(G_1,\eta,\mathbb{Q})$ be the simply connected Lie group
with Lie algebra $L_1$ equipped with the left invariant qc structure $(\eta,%
\mathbb{Q})$ defined above. Then

\begin{itemize}
\item[a)] The torsion endomorphism of the Biquard connection is zero and the normalized
qc scalar curvature is a negative constant, $S=-\frac12$.

\item[b)] The qc conformal curvature is zero, $W^{qc}=0$, and therefore $%
(G_1,\eta,\mathbb{Q})$ is locally qc conformally flat.
\end{itemize}
\end{thrm}

\begin{proof}
We compute the connection 1-forms and the horizontal Ricci forms of the Biquard connection. The Lie algebra
structure equations \eqref{ex11} together with \eqref{coneforms}, \eqref{coneform1} and \eqref{sp1curv} imply
\begin{equation}  \label{ex1conf}
\alpha_i=(\frac14-\frac{S}2)\eta_i, \qquad \rho_i(X,Y)=\frac12d\alpha_i(X,Y)=(\frac14-\frac{S}2)\omega_i(X,Y).
\end{equation}
Compare \eqref{ex1conf} with \eqref{sixtyfour} to conclude that
the torsion is zero and the  normalized qc scalar $S=-\frac12$ and
Theorem \ref{sixtyseven} completes the proof of a).

In view of Theorem~\ref{main1}, to prove b) we have to show $W^{qc}=0$. We claim $WR=0$. Indeed, since the
torsion of the Biquard connection vanishes and $S=-\frac12$, \eqref{qcwdef1} takes the form
\begin{multline}  \label{qcwdef2}
WR(X,Y,Z,V)= R(X,Y,Z,V) \\
-\frac{1}8\Big[(g\owedge g)(X,Y,Z,V)+\sum_{s=1}^3\Bigl((\omega_s\owedge%
\omega_s)(X,Y,Z,V) +4\omega_s(X,Y)\omega_s(Z,V)\Bigr) \Big].
\end{multline}
Let $A,B,C\in\Gamma(TG_1)$. Since the Biquard connection preserves the whole metric, it is connected with the
Levi-Civita connection $\nabla^g$ of the metric $g$ by the general formula
\begin{equation}  \label{lcbi}
g(\nabla_AB,C)=g(\nabla^g_AB,C)+\frac12\Big[%
g(T(A,B),C)-g(T(B,C),A)+g(T(C,A),B)\Big].
\end{equation}
The Koszul formula for a left-invariant vector fields reads
\begin{equation}  \label{lilc}
g(\nabla^g_{e_{a}}e_b,e_c)=\frac12\Big[%
g([e_a,e_b],e_c))-g([e_b,e_c],e_a)+g([e_c,e_a],e_b)\Big].
\end{equation}
Theorem~\ref{biqcon} supplies the formula
\begin{equation}  \label{torhor}
T(X,Y)=2\sum_{s=1}^3\omega_s(X,Y)\xi_s.
\end{equation}

Using \eqref{torhor}, \eqref{lilc}, \eqref{lcbi} and the structure equations %
\eqref{ex11} we found that the  non zero  Christoffel symbols for the Biquard connection (defined by
$\nabla_{e_a}e_b=\sum_c\Gamma_{ab}^ce_c$) are:

\begin{equation*}
\begin{aligned}
1&=\Gamma_{22}^1=\Gamma_{23}^4=\Gamma_{33}^1=\Gamma_{34}^2=\Gamma_{42}^3=\Gamma_{44}^1=-\Gamma_{21}^2
=-\Gamma_{24}^3=-\Gamma_{31}^3=-\Gamma_{32}^4=-\Gamma_{41}^4=-\Gamma_{43}^2,\\
\frac 12&=\Gamma_{53}^4=\Gamma_{56}^7=\Gamma_{64}^2=\Gamma_{67}^5=\Gamma_{72}^3= \Gamma_{75}^6= -\Gamma_{54}^3=
-\Gamma_{57}^6=-\Gamma_{62}^4=-\Gamma_{65}^7=-\Gamma_{73}^2=-\Gamma_{76}^5.
\end{aligned}
\end{equation*}

And the non zero coefficients of the curvature tensor are $R(e_a,e_b,e_a,e_b)=-R(e_a,e_b,e_b,e_a)=1$,
$a,b=1,\dots,4$, $a\not=b$. Now \eqref{qcwdef2} yields $\qquad WR(e_a,e_b,e_c,e_d)=R(e_a,e_b,e_a,e_b)=0$, when
there are three different indices in $a,b,c,d$. For the indices repeated in pairs we have
\begin{multline*}
WR(e_a,e_b,e_a,e_b)=R(e_a,e_b,e_a,e_b)-\frac{1}8(g\owedge %
g)(e_a,e_b,e_a,e_b)- \\
\frac{1}8\Big[\sum_{s=1}^3\Bigl((\omega_s\owedge\omega_s)(e_a,e_b,e_a,e_b)
+4\omega_s(e_a,e_b)\omega_s(e_a,e_b)\Bigr) \Big]=1-\frac18.2-\frac18.6=0
\end{multline*}
Then Theorem \ref{main1} completes the proof.
\end{proof}

\subsection{Zero torsion qc-non-flat-Example 2.}

Consider the simply connected Lie group $L_2$ with Lie algebra defined by the equations:
\begin{equation}  \label{ex2}
\begin{aligned}
& de^1 = 0, \quad de^2 = -e^{12} + e^{34},\quad de^3 = -\frac12
e^{13},\quad de^4 =-\frac12 e^{14},\\
& de^5 = 2 e^{12} + 2 e^{34} + e^{37} - e^{46} +\frac14e^{67}, \quad de^6 = 2 e^{13} - 2 e^{24} -\frac12 e^{27}
+ e^{45}
-\frac14 e^{57}, \\
& de^7 = 2 e^{14} + 2 e^{23} + \frac12 e^{26} - e^{35} +\frac14 e^{56}.
\end{aligned}
\end{equation}
A global qc structure on $L_2$ is defined by setting
\begin{equation}  \label{qc2}
\begin{aligned} &\eta_1=e^5, \quad \eta_2=e^6, \quad \eta_3=e^7, \quad
\xi_1=e_5,\quad \xi_2=e_6,\quad \xi_3=e_7,\\ &\mathbb H=span\{e^1,\dots, e^4\}, \quad \omega_1=e^{12}+e^{34},
\qquad \omega_2=e^{13}+e^{42} \qquad \omega_3=e^{14}+e^{23}, \end{aligned}
\end{equation}
It is straightforward to check from \eqref{ex2} that the triple $\{\xi_1, \xi_2, \xi_3\}$ forms the Reeb vector
fields satisfying \eqref{bi1} and therefore the Biquard connection do exists.

\begin{thrm}
\label{m22} Let $(G_2,\eta,\mathbb{Q})$ be the simply connected Lie group
with Lie algebra $L_2$ equipped with the left invariant qc structure $(\eta,%
\mathbb{Q})$ defined above. Then:

\begin{itemize}
\item[a)] The torsion endomorphism of the Biquard connection is zero and the normalized
qc scalar curvature is a negative constant, $S=-\frac14$.

\item[b)] The qc conformal curvature is not zero, $W^{qc}\not=0$
and  therefore $(G_2,\eta,\mathbb{Q})$ is not locally qc
conformally flat.
\end{itemize}
\end{thrm}

\begin{proof}
We compute the connection 1-forms and the horizontal Ricci forms of the Biquard connection. The Lie algebra
structure equations \eqref{ex2} together with \eqref{coneforms}, \eqref{coneform1} and \eqref{sp1curv} imply
\begin{equation}  \label{ex2conf}
\begin{aligned}& \alpha_1=-\frac12e^2-(\frac18+\frac{S}2)\eta_1, \quad
\alpha_2=-e^3-(\frac18+\frac{S}2)\eta_2,\quad
\alpha_3=-e^4-(\frac18+\frac{S}2)\eta_3,\\
&\rho_i(X,Y)=(\frac18-\frac{S}2)\omega_i(X,Y). \end{aligned}
\end{equation}
Compare \eqref{ex2conf} with \eqref{sixtyfour} to conclude that
the torsion is zero and the  normalized qc scalar $S=-\frac14$.
Theorem~\ref{sixtyseven} completes the proof of a).

In view of the proof of Theorem~\ref{m1}, to get b) we have to show $%
WR(e_1,e_2,e_3,e_4)=R(e_1,e_2,e_3,e_4)\not=0$.

Indeed, using \eqref{torhor}, %
\eqref{lilc}, \eqref{lcbi} and the structure equations \eqref{ex2} we found that the non zero Christoffel
symbols for the Biquard connection  are

\begin{equation*}
\begin{aligned}
1=&-\Gamma_{21}^2=\Gamma_{22}^1=\Gamma_{23}^4=-\Gamma_{24}^3=-\Gamma_{31}^3=-\Gamma_{32}^4=
\Gamma_{33}^1=\Gamma_{34}^2=-\Gamma_{41}^4=\Gamma_{42}^3=-\Gamma_{43}^2=\Gamma_{44}^1,\\
\frac12=&\Gamma_{53}^4=-\Gamma_{54}^3=\Gamma_{56}^7=-\Gamma_{57}^6=-\Gamma_{62}^4=\Gamma_{64}^2=-\Gamma_{65}^7
=\Gamma_{67}^5=\Gamma_{72}^3=-\Gamma_{73}^2=\Gamma_{75}^6=-\Gamma_{76}^5,
\end{aligned}
\end{equation*}
and then $R(e_1,e_2,e_3,e_4)=-\frac 12\not=0$. Theorem~\ref{main1} completes the proof.
\end{proof}

\subsection{Non-zero torsion qc-non-flat-Example 3}

Consider the Lie algebra defined by the equations:
\begin{equation}  \label{ex3}
\begin{aligned}
&d\tilde e^1 =\tilde e^{13} -\tilde e^{24}; \quad d\tilde e^2 = \tilde e^{14} +\tilde e^{23};\quad
d\tilde e^3 = d\tilde e^4 = 0;\\
&d\tilde e^5 = -2\tilde e^{12} - 2\tilde e^{34} - \frac12\tilde e^{17} +\frac 12\tilde e^{26}
-\tilde e^{35} -\frac18\tilde e^{67};\\
&d\tilde e^6 = -2\tilde e^{13} + 2\tilde e^{24} - \frac 12\tilde e^{36} +\frac 12\tilde
e^{47};\\
&d\tilde e^7 = -2\tilde e^{14} - 2\tilde e^{23} - \frac 12\tilde e^{37} -\frac 12\tilde e^{46}.
\end{aligned}
\end{equation}

Let $L_3$ be the Lie algebra isomorphic to \eqref{ex3} described by
\begin{equation}  \label{ex31}
\begin{aligned}
& de^1=-\frac32e^{13}+\frac32e^{24}-\frac34e^{25}+\frac14e^{36}-\frac14e^{47}+%
\frac18e^{57},\\
& de^2=-\frac32e^{14}-\frac32e^{23}+\frac34e^{15}+\frac14e^{37}+\frac14e^{46}-%
\frac18e^{56},\\
& de^3=0, \qquad
de^4=e^{12}+e^{34}+\frac12e^{17}-\frac12e^{26}+\frac14e^{67}, \\
& de^5=2e^{12}+2e^{34}+e^{17}-e^{26}+\frac12e^{67},\\ & de^6=2e^{13}+2e^{42}+e^{25},\qquad
de^7=2e^{14}+2e^{23}-e^{15}. \end{aligned}
\end{equation}
and $e_l, 1 \leq l\leq 7$ be the left invariant vector field dual to the 1-forms ${e^i, 1 \leq i\leq 7}$,
respectively. We define a global qc structure on $L_3$ by setting
\begin{equation}  \label{qc3}
\begin{aligned} &\eta_1=e^5, \quad \eta_2=e^6, \quad \eta_3=e^7, \quad
\xi_1=e_5,\quad \xi_2=e_6,\quad \xi_3=e_7,\\ &\mathbb H=span\{e^1,\dots, e^4\}, \quad \omega_1=e^{12}+e^{34},
\qquad \omega_2=e^{13}+e^{42} \qquad \omega_3=e^{14}+e^{23}, \end{aligned}
\end{equation}
It is straightforward to check from \eqref{ex3} that the triple $\{\xi_1, \xi_2, \xi_3\}$ forms the Reeb vector
fields satisfying \eqref{bi1} and therefore the Biquard connection do exists.

\begin{thrm}
\label{m3} Let $(G_3,\eta,\mathbb{Q})$ be the simply connected Lie group
with Lie algebra $L_3$ equipped with the left invariant qc structure $(\eta,%
\mathbb{Q})$ defined by \eqref{qc3}. Then

\begin{itemize}
\item[a)] The torsion endomorphism of the Biquard connection is not zero and therefore $%
(G_3,\eta,\mathbb{Q})$ is not locally qc homothetic to a 3-Sasaki manifold. The normalized qc scalar curvature
is negative, $S=-1$.

\item[b)] The qc conformal curvature is not zero, $W^{qc}\not=0$, and
therefore $(G_3,\eta,\mathbb{Q})$ is not locally qc conformally flat.
\end{itemize}
\end{thrm}

\begin{proof}
It is clear from \eqref{ex31} that the vertical distribution spaned by $%
\{\xi_1,\xi_2,\xi_3\}$ is not integrable. Consequently, the torsion of the Biquard connection is not zero due
to \cite[Theorem~3.1]{IMV} which proves the first part of a).

To prove $S=-1$  we compute the torsion. The Lie algebra structure equations \eqref{ex31} together with
\eqref{coneforms}, \eqref{coneform1} imply
\begin{equation}  \label{ex3conf}
\alpha_1=(\frac14-\frac{S}2)\eta_1, \qquad \alpha_2=-e^1-(\frac14+\frac{S}%
2)\eta_2,\qquad \alpha_3=-e^2-(\frac14+\frac{S}2)\eta_3.
\end{equation}
Now, \eqref{ex3conf}, \eqref{ex31} and \eqref{sp1curv} yield
\begin{equation}  \label{rtor}
\begin{aligned}
\rho_1(X,Y)=\frac12\Big[(\frac12-S)(e^{12}+e^{34})+e^{12}\Big](X,Y)=%
\frac14(e^{12}-e^{34})(X,Y)+ \frac12(1-S)\omega_1(X,Y),\\
\rho_2(X,Y)=\frac12\Big[\frac32(e^{13}-e^{24})(X,Y)-(%
\frac12+S)(e^{13}-e^{24})(X,Y)\Big]= +\frac12(1-S)\omega_2(X,Y),\\
\rho_3(X,Y)=\frac12\Big[\frac32(e^{14}+e^{23})(X,Y)-(%
\frac12+S)(e^{14}+e^{23})(X,Y)\Big]= +\frac12(1-S)\omega_3(X,Y) \end{aligned}
\end{equation}
Compare \eqref{rtor} with \eqref{sixtyfour} to conclude
\begin{equation}  \label{tr3}
\begin{aligned} T^0(X,I_1Y)-T^0(I_1X,Y)=\frac12(e^{12}-e^{34})(X,Y), \qquad
S=-1\\ T^0(X,I_2Y)-T^0(I_2X,Y)=0, \quad T^0(X,I_3Y)-T^0(I_3X,Y)=0.
\end{aligned}
\end{equation}
To prove b) we compute  the tensor $WR$.
Denote $\psi=-\frac14(e^{12}-e^{34})$ and compare \eqref{tr3} with %
\eqref{propt} and \eqref{need} to obtain
\begin{equation}  \label{tor3}
T^0(X,Y)=\psi(X,I_1Y), \qquad g(T(\xi_s,X),Y)=-\frac14(\psi(I_sX,I_1Y)+\psi(X,I_1I_sY)).
\end{equation}
Using $U=0$ and  \eqref{propt} we conclude  from %
\eqref{qcwdef1} that $
WR(e_1,e_2,e_3,e_4)=R(e_1,e_2,e_3,e_4)$ since other terms on the right hand
side of \eqref{qcwdef1} vanish on the quadruple $%
\{e_1,e_2=-I_1e_1,e_3=-I_2e_1,e_4=-I_3e_1\}$.

We calculate $R(e_1,e_2,e_3,e_4)$ using \eqref{lcbi}, \eqref{lilc}, %
\eqref{torhor}, \eqref{ex31} and \eqref{tor3}. We have
\begin{equation*}
\begin{aligned}
&\frac32=\Gamma_{13}^1=-\Gamma_{11}^3=-\Gamma_{22}^3=\Gamma_{23}^2,\quad\quad
\frac12=-\Gamma_{12}^4=\Gamma_{14}^2=\Gamma_{21}^4=-\Gamma_{24}^1,\\
&\frac34=\Gamma_{51}^2=-\Gamma_{52}^1=\Gamma_{56}^7=-\Gamma_{57}^6,\quad\quad
\frac18=-\Gamma_{61}^3=\Gamma_{63}^1=-\Gamma_{72}^3=\Gamma_{73}^2,\\
&\frac14=-\Gamma_{65}^7=\Gamma_{67}^5=\Gamma_{75}^6=-\Gamma_{76}^5,\quad\quad
\frac38=-\Gamma_{62}^4=\Gamma_{64}^2=\Gamma_{71}^4=-\Gamma_{74}^1,\\
&1=\Gamma_{15}^7=-\Gamma_{17}^5=-\Gamma_{25}^6=\Gamma_{26}^5=-\Gamma_{41}^2=
\Gamma_{42}^1=\Gamma_{43}^4=-\Gamma_{44}^3.
\end{aligned}
\end{equation*}
This gives $WR(e_1,e_2,e_3,e_4)=R(e_1,e_2,e_3,e_4)=-\frac12\not=0$ and Theorem \ref{main1} completes the proof.
\end{proof}

\section{$Sp(n)Sp(1)$-hypo structures and hypersurfaces in quaternionic K\"ahler manifolds}
\label{qsph}

Guided by the Examples 1-3, we relax the definition of a qc structure dropping the ``contact condition''
$d\eta_{s_{|H}}=2\omega_s$ and
come to an $Sp(n)Sp(1)$ structure (almost 3-contact structure see \cite{Kuo}%
). The purpose is to get a structure which possibly may induce an explicit quaternionic K\"ahler metric on a
product with a real line.

\begin{dfn}\label{d:sp(n)sp(1)structure}
An $Sp(n)Sp(1)$ structure on a $(4n+3)$-dimensional Riemannian
manifold $(M,g)$ is a codimension three distribution $H$
satisfying
\begin{enumerate}
\item[i)] { $H$ has an $Sp(n)Sp(1)$ structure, that is it is
equipped with a Riemannian metric $g$ and a rank-three bundle
$\mathbb Q$ consisting of (1,1)-tensors on $H$ locally generated
by three almost complex structures $I_1,I_2,I_3$ on $H$ satisfying
the identities of the imaginary unit quaternions,
$I_1I_2=-I_2I_1=I_3, \quad I_1I_2I_3=-id_{|_H}$ which are
hermitian compatible with the metric $g(I_s.,I_s.)=g(.,.), s =
1,2,3$, i.e. $H$ has an almost quaternion hermitian structure.}
\item[ii)] $H$ is locally given as the kernel of a 1-form
$\eta=(\eta_1,\eta_2,\eta_3)$ with values in $\mathbb{R}^3$.
\end{enumerate}
The local fundamental 2-forms are defined on $H$ as usual by $\omega_s(X,Y)=g(I_sX,Y)$.
\end{dfn}

\begin{dfn} We define a global $Sp(n)Sp(1)$-invariant
4-form of an $Sp(n)Sp(1)$ structure $(M,g,\mathbb Q)$ on a $(4n+3)$-dimensional manifold $M$ by the expression
\begin{equation}\label{res4}
\Omega=\omega_1^2+\omega_2^2+\omega_3^2+2\omega_1\wedge\eta_2\wedge\eta_3+
2\omega_2\wedge\eta_3\wedge\eta_1+2\omega_3\wedge\eta_1\wedge\eta_2.
\end{equation}
\end{dfn}

Let $M^{4n+4}$ be a $(4n+4)$- dimensional manifold equipped with an $Sp(n+1)Sp(1)$%
 structure, i.e.\\ $(M^{4n+4},g,J_1,J_2,J_3)$ is
{ an almost quaternion hermitian} manifold with local K\"ahler
forms $F_i=g(J_i.,.)$. The  fundamental 4-form
\begin{equation}\label{4-form}\Phi=F_1\wedge F_1+F_2\wedge
F_2+F_3\wedge F_3
\end{equation}
is globally defined and encodes fundamental properties of the structure. If the holonomy of the Levi-Civita
connection is contained in $Sp(n+1)Sp(1)$  then the manifold is a quaternionic K\"ahler  manifold which is
consequently an Einstein manifold. Equivalent conditions are either that
\begin{equation}\label{qkcon}dF_i\in
span\{F_i,F_j,F_k\}
\end{equation}
\cite{Sw} or the fundamental 4-form $\Phi$ is parallel with respect to the Levi-Civita connection. The latter
is equivalent to the fact that the fundamental 4-form is closed ($d\Phi=0$) provided the dimension is strictly
bigger than eight (\cite{Sw,Sal}) with a counter-example in dimension eight constructed by Salamon in
\cite{Sal}.

Let $f: N^{4n+3}\longrightarrow M^{4n+4}$ be an oriented hypersurface of $M^{4n+4}$
and denote by $\mathbb{N}$ the unit normal vector field. Then an $Sp(n+1)Sp(1)$
structure on $M$ induces an $Sp(n)Sp(1)$ structure on $N^{4n+3}$ locally given by $(\eta_s,\omega_s)$  defined
by the equalities
\begin{equation}  \label{qhyp1}
\eta_s=\mathbb{N}\lrcorner F_s,\quad \omega_i=f^*F_i-\eta_j\wedge\eta_k,
\end{equation}

\noindent for any cyclic permutation $(i,j,k)$ of $(1,2,3)$. The fundamental four form $\Phi$ on $M$ restricts
to the fundamental four form $\Omega$ on $N$,
\begin{equation}\label{restr4}
\Omega=f^*\Phi=(f^*F_1)^2+(f^*F_2)^2+(f^*F_3)^2.
\end{equation}
Suppose that $(M^{4n+4},g)$ has holonomy contained in $Sp(n+1)Sp(1)$. Then $d\Phi=0$, \eqref{restr4} and
\eqref{qhyp1} imply that the $Sp(n)Sp(1)$ structure  induced on $N^{4n+3}$ satisfies the equation
\begin{equation}  \label{qrhypo}
d\Omega=0,
\end{equation}
since $d$ comutes with $f^*$, $df^*=f^*d$.
\begin{dfn}
An $Sp(n)Sp(1)$ structure $(M,g,\mathbb Q)$ on a $(4n+3)$-dimensional manifold is called
\textrm{Sp(n)Sp(1) - hypo} if its fundamental 4-form is closed, $d\Omega=0$.
\end{dfn}

Hence, any oriented hypersurface $N^{4n+3}$ of a quaternionic K\"ahler $M^{4n+4}$ is naturally endowed with an
$Sp(n)Sp(1)$-hypo structure.

Vice versa, a $(4n+3)$-manifold $N^{4n+3}$ with an $Sp(n)Sp(1)$ structure $%
(\eta_s,\omega_s)$ induces an $Sp(n+1)Sp(1)$ structure $(F_s)$ on $N^{4n+3}\times%
\mathbb{R}$ defined by
\begin{equation}  \label{qhyp0}
F_i=\omega_i+\eta_j\wedge\eta_k-\eta_i\wedge dt,
\end{equation}
where $t$ is a coordinate on $\mathbb{R}$.

Consider $Sp(n)Sp(1)$ structures $(\eta_s(t),\omega_s(t))$ on $N^{4n+3}$
depending on a real parameter $t\in\mathbb{R}$, and the corresponding $%
Sp(n+1)Sp(1)$ structures $F_s(t)$ on $N^{4n+3}\times\mathbb{R}$. We have

\begin{prop}\label{qevpro} An $Sp(n)Sp(1)$ structure $(\eta_s,\omega_s; 1\leq s \leq
3)$ on $N^{4n+3}$ can be lifted to a quaternionic K\"ahler structure $(F_s(t))$ on $N^{4n+3}\times\mathbb{R}$
defined by \eqref{qhyp0} if and only if it is an $Sp(n)Sp(1)$-hypo structure
which generates a 1-parameter family of $Sp(n)Sp(1)$-hypo structures $(\eta_s(t),%
\omega_s(t))$ satisfying the following \emph{evolution Sp(n)Sp(1)-hypo equations}
\begin{equation}  \label{qevolunk}
\partial_t\Omega(t)=d\Big[6\eta_1(t)\wedge\eta_2(t)\wedge\eta_3(t)+2\omega_1(t)\wedge\eta_1(t)+
2\omega_2(t)\wedge\eta_2(t)+2\omega_3(t)\wedge\eta_3(t)\Big],
\end{equation}
where $d$ is the exterior derivative on $N$.
\end{prop}

\begin{proof}
If we apply \eqref{qhyp0} to \eqref{4-form} and then take  the exterior derivative in the obtained equation we
see that the equality $d\Phi=0$ holds precisely when \eqref{qrhypo} and \eqref{qevolunk} are fulfilled.

It remains to show that the equations \eqref{qevolunk} imply that %
\eqref{qrhypo} hold for each $t$. Indeed, using \eqref{qevolunk}, we calculate
\begin{gather*}
\partial_td\Omega= d^2\Big[6\eta_1(t)\wedge\eta_2(t)\wedge\eta_3(t)+2\omega_1(t)\wedge\eta_1(t)+
2\omega_2(t)\wedge\eta_2(t)+2\omega_3(t)\wedge\eta_3(t)\Big]=0.
\end{gather*}
Hence, the equalities \eqref{qrhypo} are independent of $t$ and therefore valid for all $t$ since it holds in
the beginning for $t=0$.
\end{proof}

Solutions to the \eqref{qrhypo} are given in the case of 3-Sasakian manifolds in \cite{WM}. In the next section
we construct explicit examples relying on the properties of the qc structures.

In general, a question remains.

\textbf{Question 1.} Does the converse of Proposition~\ref{qevpro} hold?, i.e. is it true that any
$Sp(n)Sp(1)$-hypo structure on $N^{4n+3}$ can be lifted to a quaternionic K\"ahler structure on
$N^{4n+3}\times\mathbb{R}$?

\subsection{$Sp(n)$-hypo structures and hypersurfaces in hyper
K\"ahler manifolds} \label{sph}

Suppose that $M^{4n+4}$ has holonomy contained in $Sp(n+1)$, that is the $%
Sp(n+1)Sp(1)$ structure $(F_s)$ is globally defined  and integrable (i.e. hyper-K\"ahler) or, equivalently due
to Hitchin \cite{Hit},
\begin{equation}  \label{clh}
dF_s=0.
\end{equation}
Then, \eqref{clh} and \eqref{qhyp1} imply that the $Sp(n)$ structure $%
(\eta_s,\omega_s)$ induced on $N^{4n+3}$ satisfies the equations
\begin{equation}  \label{rhypo}
d(\omega_i+\eta_j\wedge\eta_k)=0,
\end{equation}
since $d$ commutes with $f^*$, $df^*=f^*d$.

\begin{dfn}
An $Sp(n)$ structure determined by $(\eta_s,\omega_s)$ on a $(4n+3)$-dimensional manifold is called
\textrm{Sp(n)-hypo} if it satisfies the equations \eqref{rhypo}
\end{dfn}

Hence, any oriented hypersurface $N^{4n+3}$ of a hyper K\"ahler $M^{4n+4}$ is naturally endowed with an
$Sp(n)$-hypo structure.

Vice versa, a $(4n+3)$-manifold $N^{4n+3}$ with an $Sp(n)$ structure $%
(\eta_s,\omega_s)$ induces an $Sp(n+1)$ structure $(F_s)$ on $N^{4n+3}\times%
\mathbb{R}$ defined by \eqref{qhyp0}.

Consider $Sp(n)$ structures $(\eta_s(t),\omega_s(t))$ on $N^{4n+3}$
depending on a real parameter $t\in\mathbb{R}$, and the corresponding $%
Sp(n+1)$ structures $F_s(t)$ on $N^{4n+3}\times\mathbb{R}$. We have

\begin{prop}
\label{evpro} An $Sp(n)$ structure $(\eta_s,\omega_s; 1\leq s \leq 3)$ on $N^{4n+3}$ can be lifted to a hyper
K\"ahler structure $(F_s(t))$ on $N^{4n+3}\times\mathbb{R}$ defined by \eqref{qhyp0} if and only if it is an
$Sp(n)$-hypo structure which
generates an 1-parameter family of $Sp(n)$ structures $(\eta_s(t),%
\omega_s(t))$ satisfying the following \emph{evolution $Sp(n)$-hypo equations}
\begin{equation}  \label{evolunk}
\partial_t(\omega_i+\eta_j\wedge\eta_k)=d\eta_i.
\end{equation}
\end{prop}

\begin{proof}
Taking the exterior derivatives in \eqref{qhyp0} shows that the equalities $dF_s=0$ hold precisely when
\eqref{rhypo} and \eqref{evolunk} are fulfilled.

It remains to show that the equations \eqref{evolunk} imply that %
\eqref{rhypo} hold for each $t$. Indeed, using \eqref{evolunk}, we calculate
\begin{gather*}
\partial_t\Big[d(\omega_i+\eta_j\wedge\eta_k)\Big]= d^2\eta_i=0.
\end{gather*}
Hence, the equalities \eqref{rhypo} are independent of $t$ and therefore valid for all $t$ since it holds in
the beginning for $t=0$.
\end{proof}

It is known, \cite{BGN}, that the cone over a 3-Sasaki
manifold is hyper-K\"ahler, i.e., there is a solution to %
\eqref{evolunk}. Indeed, for a 3-Sasaki manifold we have \cite{IMV} $S=2,
d\eta_i(\xi_j,\xi_k)=2, \alpha_s=-2\eta_s$ and  the structure equations %
\eqref{streq} of a 3-Sasaki manifold become  \eqref{streq3}.
A solution to  \eqref{evolunk} is given by  $F_i(t)=t^2\omega_i+t^2\eta_j%
\wedge\eta_k-t\eta_i\wedge dt$.

In general, a question remains.

\textbf{Question 2.} Does the converse of Proposition~\ref{evpro} hold?, i.e. is it true that any  $Sp(n)$-hypo
structure on $N^{4n+3}$ can be lifted to a hyper K\"ahler structure on $N^{4n+3}\times\mathbb{R}$?

\begin{rmrk}
Question 2 is an embedding problem analogous to the  (hypo) $SU(n)$ embedding problem solved in
\cite{ConS,Con}. Here, we consider hyper K\"ahler manifolds instead of Calabi-Yau manifolds. Since $Sp(n)$ is
contained in $SU(2n)$, it follows that an $Sp(n)$ structure $(\omega_i,\eta_i)$ on  a $(4n+3)$-dimensional
manifold induces an $SU(2n)$ structure $(\eta_1,F,\Omega)$ where the $2$-form $F$
and the complex $(2n+1)$-form $\Omega$ are defined by $F=\omega_1+\eta_2\wedge\eta_3,\qquad
\Omega=(\omega_2+\sqrt{-1}\omega_3)^n\wedge(\eta_2+\sqrt{-1}\eta_3)$. Direct computations show that the
$Sp(n)$-hypo conditions \eqref{rhypo} yield  the $SU(2n)$-hypo conditions $dF=0,\quad d(\eta_1\wedge\Omega)=0$
which, in the real analytic case, imply an embedding into a Calabi-Yau manifold \cite{ConS,Con}. Hence, it
follows that any real analytic $(4n+3)$-manifold with an $Sp(n)$-hypo structure can be embedded in a Calabi-Yau
manifold. However, it is not clear whether this Calabi-Yau structure is a hyper K\"ahler.

A proof of the embedding property could be achieved following the considerations in the recent paper by Diego
Conti \cite{Con}. Consider $Sp(n)$ as a subgroup of $SO(4n)$, one has to show the existence of an ${\mathbb
I}^{Sp(n)}$-ordinary flag in the sense of \cite{Con}.
\end{rmrk}

\section{Examples of quaternionic K\"ahler structures}
In this section we suppose that $M$ is a Riemannian manifold of dimension $4n+3$ equipped with an $Sp(n)Sp(1)$
structure as in Definition \ref{d:sp(n)sp(1)structure}. We shall denote with $g_H$ the metric on the horizontal
distribution $H$. In addition, we assume that for some constant $\tau$ the following structure equations hold
\begin{equation}\label{e:structure equation for qK}
d\eta_i = 2\omega_i + 2\tau \eta_j\wedge\eta_k,
\end{equation}
for any cyclic permutation $(i,j,k)$ of $(1,2,3)$. Examples of such manifolds are provided by the following
quaternionic contact manifolds: i) the quaternionic Heisenberg group,
where $\tau=0$;
 ii) any 3-Sasakian
manifold, where $\tau=1$ (see \cite{IV1}
where it is proved that these structure equations
characterize the 3-Sasakian quaternionic contact manifolds); and
iii) the zero torsion qc-flat group $G_1$ defined in
Theorem~\ref{m1} with the structure equations described in
\eqref{ex11}, where $\tau=-1/4$. Actually, this is the only Lie
group satisfying the structure equation \eqref{e:structure
equation for qK} for some (necessarily) negative constant $\tau$.
We prefer to include the parameter $\tau$ since it describes
qc structures homothetic to each other. In particular, for
$\tau<0$ (resp. $\tau>0$), the qc homothety $\eta_i \mapsto \,
-{2\tau}\eta_j$ (resp. $\eta_i \mapsto \, {\tau}\eta_j$) brings
the  qc structure \eqref{qc1} on the lie algebra \eqref{ex11}
(resp. a 3-Sasakain structure) to one satisfying
\eqref{e:structure equation for qK}. On the other hand, this one
parameter family of homothetic to each other qc structures lead to
different special holonomy metrics, which we construct next, when
we take the product with a real line.


\begin{thrm}\label{qkmetric}
Let $M$ be a smooth manifold of dimension $4n+3$ equipped with an $Sp(n)Sp(1)$ structure such that, for some
constant $\tau$, the structure equations \eqref{e:structure equation for qK} hold for any cyclic permutation
$(i,j,k)$ of $(1,2,3)$. For any constant $a$, the manifold $M\times\mathbb{R}$ has a quaternionic K\"ahler
structure given by the following metric and fundamental 4-form
\begin{equation}\label{e:general qK metric}
\begin{aligned}
&g=ug_H+(\tau u+au^2)(\eta_1^2+\eta_2^2+\eta_3^2)+\frac {1}{4(\tau u+au^2)}(du)^2,\quad \tau u+au^2>0,\\
&\Phi=F_1\wedge F_1 +F_2\wedge F_2 + F_3\wedge F_3,
\end{aligned}
\end{equation}
where locally
\begin{equation}\label{e:local qK forms}
\begin{aligned}
F_i(u)=u\omega_i+(a u^2+\tau u)\,\eta_j\wedge\eta_k-\frac 12 \eta_i\wedge du.
\end{aligned}
\end{equation}
\end{thrm}

\begin{proof}
Let $h$ and $f$ be some functions of the unknown $t$  and $F_i(t)=f(t)\omega_i+h^2(t)\eta_j\wedge
\eta_k-h(t)\eta_i\wedge dt$ and $\Phi$ be as in \eqref{e:general qK metric}. A direct calculation shows that
($\Sigma_{(ijk)}$ means the cyclic sum)
\begin{equation*}
d\Phi=\Sigma_{(ijk)}\left [ \left ( (f^2)'-4fh\right)\omega_i\wedge\omega_i\wedge dt+\left ( 2\left (fh^2
\right )' +4\tau fh-12h^3 \right )\omega_i\wedge\eta_j\wedge\eta_k\wedge dt \right ].
\end{equation*}
Thus, if we take $h=\frac 12 f'$ we come to $$d\Phi=f'\Sigma_{(ijk)}(-f'^2+f f''+2\tau
f)\omega_i\wedge\eta_j\wedge\eta_k\wedge dt,$$ which shows that $\Phi$ is closed when
\begin{equation}\label{solqk7}
ff''-f'^2+2\tau f=0, \qquad h=\frac 12 f'.
\end{equation}
With the help of the substitution $v=-\ln f$ we see that $\left ( \frac {dv}{dt}\right )^2=4\tau e^v +4a$ for
any constant $a$. This shows that $\left(  \frac {dt}{df}  \right)^2 =\left( \frac {dt}{dv} \right )^2 \left(
\frac {dv}{df}\right )^2=\frac {1}{4(\tau f + af^2)} > 0$ and $h^2=\tau f + af^2$. Renaming $f$ to $u$ gives
the quaternionic structure in the local form \eqref{e:local qK forms} and the metric in \eqref{e:general qK
metric}. In order to see that $<F_1, F_2, F_3>$ is a differential ideal we need to compute the differentials
$dF_i$. A small calculation shows
\begin{equation}\label{e:closed ideal}
dF_i=\frac {f}{f'}(ff''-f'^2+2\tau f)\eta_j\wedge\eta_k\wedge dt \qquad mod \quad <F_1, F_2, F_3>,
\end{equation}
i.e. \eqref{qkcon} hold. This proves that the defined structure is quaternionic K\"ahler taking into account
the differential equation  \eqref{solqk7} satisfied by $f$,  which completes the proof.
\end{proof}

With the help of the above theorem we obtain the following one parameter families of quaternionic K\"ahler
structures.

i) \textit{Quaternionic K\"ahler metrics from the quaternionic Heisenberg group, $\tau=0$}. Consider the
$(4n+3)$-dimensional quaternionic Heisenberg group $\mathbb{G}^n$, viewed as a quaternionic contact structure.
The metric
\begin{equation}\label{e:aqKahler}
g=e^{at}\left ((e^1)^2+\dots +(e^{4n})^2 \right )+\frac {a^2}{4}e^{2at}\left ( (\eta_1)^2+(\eta_2)^2+(\eta_3)^2
 \right )+dt^2
\end{equation}
is a complete quaternionic K\"ahler metric in dimensions $4n+4$ with $n\geq 1$.
 The Einstein constant is negative equal to $-16na^2$. This complete Einstein metric has  been found in dimension
 eight as an Einstein metric on a $T^3$ bundle over $T^4$ in \cite[equation (148)]{GLPS}.

ii) \textit{Quaternionic K\"ahler metrics from a 3-Sasakian structure, $\tau=1$}. The metric
\begin{equation*}
g= u g_H +\frac {u+au^2}{4}\left ( (\eta_1)^2+(\eta_2)^2+(\eta_3)^2 \right )+ \frac {1}{4(u+au^2)}du^2
\end{equation*}
is a quaternionic K\"ahler, and in the case of $a=0$ is the hyper-K\"ahler cone over the 3-Sasakian manifold.
These metrics have been found earlier in \cite[Theorem 5.2]{WM}.

\subsection{Explicit non quaternionic K\"ahler structures with closed four form in dimension 8}
As it is well known \cite{Sw} in dimension $4n$, $n>2$, the
condition that the fundamental 4-form is closed is equivalent to
the fundamental 4-form being parallel which is not true in
dimension eight. Salamon constructed in \cite{Sal} a compact
example of
{an almost quaternion} hermitian manifold with closed fundamental
four form which is not Einstein, and therefore it is not
quaternionic K\"ahler.   We give below explicit complete
non-compact examples of that kind inspired by the following


\begin{rmrk}\label{r:dim seven QK}
In dimension seven, due to the relations $\omega_i\wedge\omega_j=0$, $i\not= j$,  a more general
evolution than the one considered in the proof of Theorem~\ref{qkmetric} can be handled. We consider the evolution
\begin{equation}\label{genevolution}
\omega_s(t)=f(t)\omega_s, \qquad \eta_s(t)=f_s(t)\eta_s, \qquad s=1,2,3,
\end{equation}
where $f,f_1,f_2,f_3$ are smooth function of $t$. Using the structure equations \eqref{e:structure equation for
qK} one easily obtain that the equation $d\Omega=0$ is satisfied and  \eqref{qevolunk} is equivalent to the
system
\begin{equation}\label{erealqk}
\begin{aligned}
3f'-2(f_1+f_2+f_3)=0,\\
(ff_2f_3)'-2 \tau f(f_1-f_2-f_3)-6 f_1f_2f_3=0,\\
(ff_1f_3)'-2\tau f(-f_1+f_2-f_3)-6 f_1f_2f_3=0,\\
(ff_1f_2)'-2\tau f(-f_1-f_2+f_3) -6 f_1f_2f_3=0.
\end{aligned}
\end{equation}
On the other hand, $<F_1, F_2, F_3>$ is a differential ideal  iff the following system holds
\begin{equation}\label{e:closed ideal qk}
f(f_if_j)'-f'f_if_j+2f_1f_2f_3-2f_if_j(f_i+f_j)+2\tau f f_if_j-2\tau f f_k=0,
\end{equation}
for any cyclic permutation $(i,j,k)$ of $(1,2,3)$. This claim follows from the fact that working $ mod <F_1,
F_2, F_3>$we have
\[
dF_i=\frac {1}{f}\left (f(f_if_j)'-f'f_if_j+2f_1f_2f_3-2f_i^2f_j-2f_if_j^2+2\tau ff_if_j-2\tau f f_k \right
)\eta_j\wedge\eta_k\wedge dt.
\]
 Taking
$f_1=f_2=f_3=h$ in \eqref{erealqk} we come, correspondingly, to the case considered in Theorem \ref{qkmetric}.
\end{rmrk}

We integrate the system \eqref{erealqk} completely when $\tau=0$. This is achieved by introducing the new variable
$du=f_1f_2f_3dt$, which allow to determine $ff_if_j=6(u+a_k)$, where $a_k$ is a constant and $(i,j,k)$ is a
permutation of $(1,2,3)$. Thus $f_i=\frac {f}{6(u+a_i)} \frac {du}{dt}$, $j=1,\, 2,\, 3$. With the help of
these three equations and the first equation of \eqref{erealqk} we come to $\frac {9}{f}\frac {df}{dt}=\frac
{du}{dt}\left (\frac {1}{u+a_1} + \frac {1}{u+a_2} + \frac {1}{u+a_3}\right )$, hence
$f^9=C^9(u+a_1)(u+a_2)(u+a_3)$ for some constant $C$. Now, the equations $ff_if_j=6(u+a_k)$ yield
$$f_i=\sqrt{\frac 6C} \left ( \frac {(u+a_j)^4(u+a_k)^4}{(u+a_i)^5} \right )^{1/9}$$ and then the definition of
$u$  shows $$dt=\left ( C/6 \right )^{3/2}\frac {du}{\left ((u+a_1)(u+a_2)(u+a_3)\right )^{1/3}}.$$ If we
impose also the system \eqref{e:closed ideal qk}, in which we substitute $2(f_i+f_j)=3f'-2f_k$ and
$f(f_if_j)'=6f_1f_2f_3-\frac 23 \left ( f_1+f_2+f_3 \right )f_if_j$ (using the equations of \eqref{erealqk} and
$\tau=0$), we see that
\[
dF_i=\frac {10f_jf_k}{3f}\left ( 2f_k-f_i-f_j \right )\eta_j\wedge\eta_k\wedge dt\qquad mod \ <F_i, F_j, F_k>.
\]
Thus, $d\Phi=0$ and  $<F_1, F_2, F_3>$ is  a differential ideal if and only if $f_1=f_2=f_3$ which yield.
\begin{prop}
The  metric on the product of the seven dimensional quaternionic
Heisenberg group with the  real line defined (on $\mathbb R^8$) by
\begin{multline}\label{e:QK using H general case}
g=C \, \left ((u+a_1)(u+a_2)(u+a_3)\right )^{1/9}\, (dx_1^2+dx_2^2+dx_3^2+dx_4^2) + \\
\frac {6}{C}\left ( \frac
{(u+a_2)^8(u+a_3)^8}{(u+a_1)^{10}} \right )^{1/9}(dx_5+2x_1dx_2+x_3dx_4)^2
+\frac {6}{C}\left ( \frac
{(u+a_3)^8(u+a_1)^8}{(u+a_2)^{10}} \right )^{1/9}(dx_6+2x_1dx_3+x_4dx_2)^2\\+
\frac {6}{C}\left ( \frac
{(u+a_1)^8(u+a_2)^8}{(u+a_3)^{10}} \right )^{1/9}(dx_7+2x_1dx_4+x_2dx_3)^2
+ \left (\frac {C}{6}\right )^3\frac {du^2}{\left ((u+a_1)(u+a_2)(u+a_3)\right )^{2/3}},
\end{multline}
where $a_1, a_2$ and $a_3$ are three constants not all of them equal to each other,
supports { an almost quaternion} hermitian structure which has closed fundamental form, but
is not quaternionic K\"ahler.
\end{prop}
\begin{rmrk}
Using a suitable computer program one can check the the metrics
\eqref{e:QK using H general case}  are Einstein exactly when
$f_1=f_2=f_3$, i.e. when are quaternionic K\"ahler.
\end{rmrk}

We note that one of the arbitrary constants in \eqref{e:QK using H
general case} is unnecessary since a translation of the  unknown
$u$ does not change the metric.

Let us also remark that the quaternionic K\"ahler metric \eqref{e:aqKahler} is obtained from the general family
\eqref{e:QK using H general case} by taking $\frac {6}{C^3}=\frac {a^2}{4}$ and $v=e^{at}=Cu^{1/3}$ when the
constants are the same $a_1=a_2=a_3$ and we use $u+a_1$ as a variable, which is denoted also by $u$.

If one takes a solution of the system \eqref{e:closed ideal qk} which does not satisfy the system
\eqref{erealqk}, one could obtain a non quaternionic Kahler manifold with { an almost quaternion hermitian} structure
such that $<F_1, F_2, F_3>$ is a differential ideal and
{the fundamental four form is non-parallel (see also the
paragraph after \cite[Corollary 2.4]{Macia}).}

\subsection{New quaternionic K\"ahler metrics from the zero-torsion qc-flat qc structure on $G_1$}
Here we consider the Lie group defined by the structure equations
\eqref{ex11}, which can be described  in local coordinates
$\{t,x,y,z,x_5,x_6,x_7\}$ as follows
\begin{equation}\label{loccoord}
\begin{aligned}
& e^1=-dt,\\
& e^2={\frac 12}\,{x_6}\,dx+{\frac 12}\,{x_5}\cos x\,dy
 +({\frac 12}\,{x_6}\cos y+{\frac 12}{x_5}\sin y\sin x)\,dz-{\frac 12}\,{x_7}\,dt+{\frac 12}\,d{x_7},\\
& e^3=-{\frac 12}\,{x_7}\,dx+{\frac 12}\,{x_5}\sin x\,dy
+(-{\frac 12}\,{x_7}\cos y-{\frac 12}{x_5}\sin y\,\cos x)\,dz-{\frac 12}\,{x_6}\,dt+{\frac 12}\,d{x_6},\\
& e^4=(-{\frac 12}{x_7}\cos x\,-{\frac 12}{x_6}\sin x\,)\,dy
-{\frac 12}\sin y\,(-{x_6}\cos x+{x_7}\sin x)\,dz-{\frac 12}\,{x_5}\,dt+{\frac 12}\,d{x_5},\\
& \eta_1=e^5=-{x_6}\,dx+(-{x_5}\cos x-2\sin x)\,dy\\ &\hskip1.7truein
+(-{x_6}\cos y-\sin y\sin x\,{x_5}+2\sin y\cos x)\,dz+{x_7}\,dt-d{x_7},\\
& \eta_2=e^6={x_7}\,dx+(2\cos x-{x_5}\sin x)\,dy\\ &\hskip1.7truein
 +({x_7}\cos y+2\sin y\sin x+{x_5}\sin y\,\cos x)\,dz+{x_6}\,dt-d{x_6},\\
& \eta_3=e^7=-2\,dx+(\cos x\,{x_7}+{x_6}\sin x)\,dy\\ &\hskip1.7truein +(-2\cos y+{x_7}\sin y\sin x-{x_6}\sin
y\,\cos x)\,dz+{x_5}\,dt-d{x_5}.
\end{aligned}
\end{equation}
In this case $\tau=-\frac14$ in \eqref{e:general qK metric}, and
the corresponding quaternionic  K\"ahler metric on $G_1$ is (using
$a/4$ as a constant)
\begin{equation}\label{qknew8}
g= u\left ((e^1)^2+(e^2)^2+(e^3)^2+(e^4)^2 \right )+\frac {au^2-u}{4}\left ( (\eta_1)^2+(\eta_2)^2+(\eta_3)^2
\right )+ \frac {1}{au^2-u}du^2,
\end{equation}
for $au^2-u>0$. The Ricci tensor is given by $Ric=-4ag$.

The metric \eqref{qknew8} seems to be a new explicit quaternionic K\"ahler metric. In local
coordinates $\{v^1=t,v^2=x,v^3=y,v^4=z,v^5=x_5,v^6=x_6,v^7=x_7,v^8=u\}$ the metric
has the expression written in Appendix 1.

\section{$Sp(1)Sp(1)$ structures and $Spin(7)$-holonomy metrics}
{ An} $Sp(1)Sp(1)$ structure on a seven dimensional manifold
$M^7$, { whose  2-forms $\omega_i$
and $1$-forms $\eta_j$ are
 globally defined,}
induces a $G_2$-form $\phi$ given by
\begin{equation}\label{g2}
\phi=2\omega_1\wedge\eta_1+2\omega_2\wedge\eta_2-2\omega_3\wedge\eta_3+2\eta_1\wedge\eta_2\wedge\eta_3.
\end{equation}
The Hodge dual $*\phi$ is
\begin{equation}\label{starg2}
*\phi=-(\omega_1\wedge\omega_1+2\omega_1\wedge\eta_2\wedge\eta_3+
2\omega_2\wedge\eta_3\wedge\eta_1-2\omega_3\wedge\eta_1\wedge\eta_2).
\end{equation}
Consider the $Spin(7)$-form $\Psi$ on $M^7\times\mathbb R$ defined by \cite{BH}
\begin{equation}\label{spin7}
\Psi=F_1\wedge F_1+F_2\wedge F_2-F_3\wedge F_3=-*\phi-\phi\wedge dt,
\end{equation}
where the 2-forms $F_1,F_2,F_3$ are given  by \eqref{qhyp0}.

Following Hitchin, \cite{Hitt}, the $Spin(7)$-form $\Psi$ is
closed if and only if the $G_2$ structure is cocalibrated,
$d*\phi=0$, and the Hitchin flow equations
$\partial_t(*\phi)=-d\phi$ are satisfied, i.e.
\begin{equation}\label{evolspin}
d(*\phi)=0, \qquad \partial_t(*\phi)=-d\phi.
\end{equation}

\begin{thrm}\label{spin7metric}
Let $M$ be a smooth seven dimensional manifold equipped with an $Sp(1)Sp(1)$ structure
such that the $3$-form $\phi$ determined by \eqref{g2}
is globally defined and, for some constant $\tau\not=0$,
 the structure equations \eqref{e:structure equation for qK} hold for any cyclic
permutation $(i,j,k)$ of $(1,2,3)$. For any constant $a$, the manifold $M\times I$,
where $I\subset\mathbb{R}$, has a parallel
$Spin(7)$ structure given by the following metric and fundamental 4-form
\begin{equation}\label{e:general Spin(7) metric}
\begin{aligned}
&g=ug_H+\frac{\tau u^{5/3}-a}{5u^{2/3}}\left ( (\eta_1)^2+(\eta_2)^2+(\eta_3)^2 \right )+ \frac
{5u^{2/3}}{36 (\tau u^{5/3}-a)}du^2,\\
&\psi=F_1\wedge F_1 + F_2\wedge F_2 - F_3\wedge F_3,
\end{aligned}
\end{equation}
where locally
\begin{equation}\label{e:local spin(7) forms}
\begin{aligned}
F_i(u)=u\omega_i-\epsilon_i\frac {\tau u^{5/3}-a}{5u^{2/3}}\,\eta_j\wedge\eta_k-\epsilon_i \frac
16\,\eta_i\wedge du.
\end{aligned}
\end{equation}
where $\epsilon_1=\epsilon_2=1$  and $\epsilon_3=-1$.
\end{thrm}

\begin{proof}
We evolve the structure as in \eqref{genevolution}.
Using the structure equations \eqref{e:structure equation
for qK} one easily obtains that the equation $d(*\phi)=0$ is satisfied, and the second equation of the system
\eqref{evolspin} is equivalent to the system
\begin{equation}\label{ereal7}
\begin{aligned}
f'-2(f_1+f_2-f_3)=0,\\
(ff_2f_3)'-2 \tau f(f_1-f_2+f_3)-2 f_1f_2f_3=0,\\
(ff_1f_3)'-2\tau f(-f_1+f_2+f_3)-2 f_1f_2f_3=0,\\
(ff_1f_2)'-2\tau f(f_1+f_2+f_3)+2 f_1f_2f_3=0.
\end{aligned}
\end{equation}
Taking $f_1=f_2=-f_3$ in \eqref{ereal7} we come to the ODE system
\begin{equation}\label{sol7}
3ff''+(f')^2-18\tau f=0, \qquad f_1=f_2=-f_3=\frac16 f'.
\end{equation}
To solve this differential equation, we use  $v=f^{4/3}$ as a variable. Equation \eqref{sol7} shows  that
$\left ( \frac {dv}{dt}\right )^2= \frac{64(\tau v^{5/4}-a)}{5}$, where $a$ is a constant. Hence, $\left( \frac
{dt}{df} \right)^2 =\left( \frac {dt}{dv} \right )^2 \left( \frac {dv}{df}\right )^2=\frac {5f^{2/3}}{36(\tau
f^{5/3}-a)}$, which implies that $f_1^2=f_2^2=f_3^2=\frac{1}{36} (f')^2 =\frac{\tau f^{5/3}-a}{5f^{2/3}}$.
Recall that we also have $F_i(t)=f(t)\omega_i+f_j(t)f_k(t)\eta_j\wedge \eta_k-f_i(t)\eta_i\wedge dt$. Renaming
$f$ to $u$ gives the metric { and the $Spin(7)$ form $\psi$ in \eqref{e:general Spin(7) metric} with
\eqref{e:local spin(7) forms}}.
\end{proof}

i) \textit{$Spin(7)$-holonomy metrics from the quaternionic Heisenberg group.} Using the seven dimensional
quaternionic Heisenberg group with structure equations \eqref{4n+3heis}, taken for $n=1$,
the corresponding eight dimensional
$Spin(7)$-holonomy metric written with respect to the parameter $u=(at+b)^{1/4}$ is
\begin{equation}\label{citsp7}
\begin{aligned}
g=u^3\left ((e^1)^2+(e^2)^2+(e^3)^2+(e^4)^2 \right )+\frac {a^2}{16}u^{-2}\left (
(\eta_1)^2+(\eta_2)^2+(\eta_3)^2 \right )+ \frac {4}{a^2}u^{6}du^2.
\end{aligned}
\end{equation}
These $Spin(7)$-holonomy metrics are found in \cite[Section
4.3.1]{GLPS}.

ii)  \textit{$Spin(7)$-holonomy metrics from a 3-Sasakian manifold.} This case was investigated
in general in
\cite{Baz} and explicit solutions in particular cases are known (see \cite{Baz} and references therein). We use
again only the particular solution to \eqref{ereal7} found above. Thus, starting with a 3-Sasakian manifold
with structure equations \eqref{streq3} the resulting metric is
\begin{equation*}
g=u\left ((e^1)^2+(e^2)^2+(e^3)^2+(e^4)^2 \right )+\frac {u^{5/3}-a}{5u^{2/3}}\left (
(\eta_1)^2+(\eta_2)^2+(\eta_3)^2 \right )+ \frac {5u^{2/3}}{36(u^{5/3}-a)}du^2.
\end{equation*}
This is the (first) complete metric with holonomy $Spin(7)$ constructed by Bryant and Salamon \cite{BrS,GPP}.

\subsection{New $Spin(7)$-holonomy metrics from the quaternionic Heisenberg group}
New metrics can be obtained similarly to the derivation of
\eqref{e:QK using H general case}. Namely, we integrate the system
\eqref{ereal7} when $\tau=0$ to obtain the next family of
$Spin(7)$-holonomy metrics   which seems to be new
\begin{multline}\label{e:Spin(7) using H general case}
g=C \, \left ((u+a_1)(u+a_2)(a_3-u)\right )\, \left
((dx_1^2+dx_2^2+dx_3^2+dx_4^2 \right ) \\+ \frac {2}{C} \frac
{1}{(u+a_1)^{2}}(dx_5+2x_1dx_2+2x_3dx_4)^2 + \frac2C\frac
{1}{(u+a_2)^{2}} (dx_6+2x_1dx_3+2x_4dx_2)^2+ \\\frac2C\frac
{1}{(a_3-u)^{2}}(dx_7+2x_1dx_4+2x_2dx_3)^2
 + \frac {C^3}{8}(u+a_1)^2(u+a_2)^2(a_3-u)^2du^2.
\end{multline}
Taking $a_2=-a_3=a_1$ into \eqref{e:Spin(7) using H general case}
one gets the Spin(7)-holonomy metrics \eqref{citsp7}. Since the
coefficients of the metrics \eqref{e:Spin(7) using H general case}
are continuous with respect to the parameters, and since the
holonomy is equal to Spin(7) for $(a_2,a_3)=(a_1,-a_1)$ then the
same holds for any $(a_2,a_3)$ in a neighbourhood of $(a_1,-a_1)$.
Thus, we get a three parameter family of metrics with holonomy
equal to $Spin(7)$ which seem to be new.

\subsection{New $Spin(7)$-holonomy metrics from a zero-torsion qc-flat qc structure on $G_1$}
Consider the $7$-dimensional Lie group defined in \eqref{ex11}. From Theorem \ref{spin7metric} we obtain the
metrics
\begin{equation}\label{newspin7metric}
g=u\left ((e^1)^2+(e^2)^2+(e^3)^2+(e^4)^2 \right )+ \frac{(a-u^{5/3})}{20u^{2/3}}\left (
(\eta_1)^2+(\eta_2)^2+(\eta_3)^2 \right )+ \frac {5u^{2/3}}{9(a-u^{5/3})}du^2.
\end{equation}
These metrics have holonomy equal to $Spin(7)$. In local
coordinates
$\{v^1=t,v^2=x,v^3=y,v^4=z,v^5=x_5,v^6=x_6,v^7=x_7,v^8=u\}$ the
$Spin(7)$-holonomy metric  is written in Appendix 2.

\section{Hyper K\"ahler metrics in dimension four}

In this section we recover some of the known  Ricci-flat gravitational instantons in dimension four  applying
our method from the preceding section lifting the $sp(0)$-hypo structures on the non-Euclidean Bianchi type
groups of class A.

Let $G_3$ is a three dimensional Lie group with Lie algebra $g_3$ and $%
e^1,e^2,e^3$ be a basis of left invariant 1-forms. We consider the $Sp(1)$ structure on $g_3\times\mathbb{R}^+$
defined by the following 2-forms
\begin{equation}  \label{4inst}
\begin{aligned} F_1(t)=e^1(t)\wedge e^2(t)+e^3(t)\wedge f(t)dt,\\
F_2(t)=e^1(t)\wedge e^3(t)-e^2(t)\wedge f(t)dt,\\
F_3(t)=e^2(t)\wedge e^3(t)+e^1(t)\wedge f(t)dt, \end{aligned}
\end{equation}
where $f(t)$ is a function of $t$ and $e^i(t)$ depend on $t$. With the help of Hitchin's theorem, it is
straightforward to prove the next

\begin{prop}
The $Sp(1)$ structure $(F_1,F_2,F_3)$ is hyper K\"ahler if and only if
\begin{equation}\label{4inst1} de^{12}=de^{13}=de^{23}=0
\end{equation}
and the following evolution equations hold
\begin{equation}\label{evol4}
\frac{\partial}{\partial t}e^{ij}(t)=-f(t)de^k(t).
\end{equation}
The hyper K\"ahler metric is given by
\begin{equation}  \label{hypkel4}
g=(e^1(t))^2+(e^2(t))^2+(e^3(t))^2+f^2(t)dt^2.
\end{equation}
\end{prop}

\medskip

\noindent
{\bf The group $SU(2)$, Bianchi type IX.}
Let $G_3=SU(2)=S^3$ be described by the structure equations
\begin{equation}  \label{su2}
de^i=-e^{jk}.
\end{equation}
In terms of Euler angles the left invariant forms $e^i$ are given by
\begin{equation}\label{eulersu2}
\begin{aligned}
e^1&=&\sin\psi d\theta-\cos\psi\sin\theta d\phi,\quad
e^2&=&\cos\psi d\theta+\sin\psi\sin\theta d\phi,\quad e^3&=&d\psi
+\cos\theta d\phi.
\end{aligned}
\end{equation}
Clearly \eqref{4inst1} are satisfied. We evolve the $SU(2)$ structure as
\begin{equation}  \label{ev}
e^s(t)=f_s(t)e^s, \quad s=1,2,3,\quad {\rm (no \quad summation \quad on \quad s)}
\end{equation}
where $f_s$ are functions of $t$.

Using the structure equations \eqref{su2} we reduce the evolution equations %
\eqref{evol4} to the following system of ODEs
\begin{equation}  \label{evol4f}
\frac{d}{d t}(f_1f_2)=ff_3,\quad \frac d{dt}(f_1f_3)=ff_2, \quad \frac d{dt}(f_2f_3)=ff_1.
\end{equation}
The system \eqref{evol4f}  is equivalent to the following 'BGPP' \cite{BGPP} system
\begin{equation}\label{ah}
\frac d{dt}f_1=f\frac{f_2^2+f_3^2-f_1^2}{2f_2f_3}, \quad \frac d{dt}f_1=f\frac{f_3^2+f_1^2-f_2^2}{2f_1f_3},
\quad \frac d{dt}f_1=f\frac{f_1^2+f_2^2-f_3^2}{2f_1f_2}.
\end{equation}
The equations \eqref{ah} admit the triaxial Bianchi IX BGPP \cite{BGPP} hyper K\"ahler metrics by taking
$f=f_1f_2f_3$ and all $f_i$ different (see also \cite{GPop}) and Eguchi-Hanson \cite{EH} hyper K\"ahler metric
when two of the functions are equal.

\subsubsection{The general solution}\label{The general solution}

With the substitution $x_{i}=(f_{j}f_{k})^{2}$, the system \eqref{evol4f} becomes
\begin{equation*}
\frac{dx_{i}}{dr}=2(x_{1}x_{2}x_{3})^{1/4},
\end{equation*}%
in terms of the parameter $dr=fdt$. Hence the functions $x_{i}$ differ by a
constant, i.e, there is a function $x(r)$ such that $%
x(r)=x_{1}+a_{1}=x_{2}+a_{2}=x_{3}+a_{3}$. The equation for $x(r)$ is
\begin{equation}
\frac{d{x}}{dr}=2\left( (x-a_{1})(x-a_{2})(x-a_{3})\right) ^{1/4},\ \text{%
i.e.,\ }dr=\frac{1}{2}\frac{1}{\left( (x-a_{1})(x-a_{2})(x-a_{3})\right) ^{1/4}}dx.  \label{e:BIX1}
\end{equation}%
If we let $g(x)=\frac{1}{2}\left( (x-a_{1})(x-a_{2})(x-a_{3})\right) ^{-1/4}$ and take into account
$x_{i}=(f_{j}f_{k})^{2}$,
we see from \eqref{evol4f} that the functions $f_{i}(x)$ satisfy
\begin{equation*}
\frac{d}{dx}\left( (x-a_{i})^{1/2}\right) =g(x)f_{i}.
\end{equation*}%
Solving for $f_{i}$ we showed that the general solution of \eqref{evol4f} is
\begin{equation}\label{e:su2}
f_{i}(x)=\frac{(x-a_{j})^{1/4}(x-a_{k})^{1/4}}{(x-a_{i})^{1/4}}%
,\quad f(t)=g\left ( x (t)\right )\,x'(t), \quad
g(x)=\frac{1}{2}\left( (x-a_{1})(x-a_{2})(x-a_{3})\right) ^{-1/4},
\end{equation}
where
$a_{1},\ a_{2}$ and $a_{3}$ are constants, and $x$ is an auxiliary
independent variable (substituting any function $x=x(t)$ gives a solution of %
\eqref{evol4f} in terms of $t$).

\subsubsection{Eguchi-Hanson instantons}

A particular solution to \eqref{evol4f} is obtained by taking $x=(t/2)^4$ and $%
a_1=a_2=\frac1{16}a$, $a_3=0$, which gives
\begin{equation}\label{ehfff}
f_1=f_2=\frac t2,\quad f_3=\frac t2\sqrt{\frac{t^4-a}{t^4}}, \quad f=\sqrt{%
\frac{t^4}{t^4-a}}.
\end{equation}
This is the Eguchi-Hanson instanton \cite{EH} with the metric given by
\begin{equation*}
g=\frac{t^2}4\Big[(e^1)^2+(e^2)^2+\Big(1-\frac{a}{t^4}\Big)(e^3)^2\Big]+ %
\Big(1-\frac{a}{t^4}\Big)^{-1}(dt)^2.
\end{equation*}

\subsubsection{Triaxial Bianchi type IX BGPP metrics \protect\cite{BGPP}}
The substitution $x=t^4$, $a_1=a^4$, $a_2=b^4$ and $a_3=c^4$ gives
\begin{equation}\label{triax}
\begin{aligned}
&f_1(t)=\frac{(t^4-b^4)^\frac14(t^4-c^4)^\frac14}{(t^4-a^4)^\frac14},\qquad
f_2(t)=\frac{(t^4-a^4)^\frac14(t^4-c^4)^\frac14}{(t^4-b^4)^\frac14},\\
&f_3=\frac{(t^4-a^4)^\frac14(t^4-b^4)^\frac14}{(t^4-c^4)^\frac14},\qquad
f(t)=\frac{2t^3}{(t^4-b^4)^\frac14(t^4-c^4)^\frac14(t^4-a^4)^\frac14}.
\end{aligned}
\end{equation}

These are the triaxial Bianchi IX metrics discovered in \cite{BGPP} (see also \cite{GPop,Gh1,Gh2}), which do
not have any tri-holomorphic U(1) isometries \cite{Gh2}. In the derivation above we avoided the use of elliptic
functions.

\medskip

\noindent {\bf The group $SU(1,1)$-Bianchi type $VIII$.}
Bianchi type $VIII$ are investigated in \cite{LP1,LP2,Lor}.

Let $G_3=SU(1,1)$ be described by the structure equations
\begin{equation}  \label{su11}
de^1=-e^{23}, \quad de^2=e^{31}, \quad de^3=-e^{12}.
\end{equation}
In terms of local coordinates the left invariant forms $e^i$ are given by
\begin{equation}  \label{su11c}
\begin{aligned} e^1&=&\sinh\psi d\theta+\cosh\psi\sin\theta d\phi,
\quad e^2&=&\cosh\psi d\theta+\sinh\psi\sin\theta d\phi,\quad
e^3&=&d\psi -\cos\theta d\phi. \end{aligned}
\end{equation}
Clearly \eqref{4inst1} are satisfied. We evolve the $SU(1,1)$ structure as in \eqref{ev}. Using the structure
equations \eqref{su11c} we reduce the evolution equations \eqref{evol4} to the following system of ODEs
\begin{equation}  \label{evol4f11}
\frac{\partial}{\partial t}(f_1f_2)=ff_3,\quad \frac{\partial}{\partial t}%
(f_1f_3)=-ff_2, \quad \frac{\partial}{\partial t}(f_2f_3)=ff_1.
\end{equation}

Solutions to the above system yield corresponding hyper K\"ahler metrics \eqref{hypkel4} indicated in
\cite{BGPP}.

\medskip

\noindent
{\bf Triaxial Bianchi type VIII metrics.}  Working as in \ref{The general solution}
we obtain the following system for the functions $x_i$
\begin{equation*}
\frac {dx_3}{dr} =\frac {dx_1}{dr}= 2(x_1x_2x_3)^{1/4}, \qquad \frac {dx_2}{dr%
}= -2(x_1x_2x_3)^{1/4}.
\end{equation*}
Solving for $f_{i}$, as in the derivation \eqref{e:su2}, we find the general solution of \eqref{evol4f11} is
\begin{equation}\label{e:su2g}
\begin{aligned}
f_{1}(x)=\frac{(x-a_{3})^{1/4}(a_{2}-x)^{1/4}}{(x-a_{1})^{1/4}}, \qquad f_{2}(x)=
\frac{(x-a_{1})^{1/4}(x-a_{3})^{1/4}}{(a_{2}-x)^{1/4}},\\
f_{3}(x)=\frac{(x-a_{1})^{1/4}(a_{2}-x)^{1/4}}{(x-a_{3})^{1/4}},\quad
f(t)=g\left ( x (t)\right )\,x'(t),\quad g(x)=\frac{1}{2}\left(
(x-a_{1})(a_{2}-x)(x-a_{3})\right) ^{-1/4},
\end{aligned}
\end{equation}
where $a_{1},\ a_{2}$ and $a_{3}$ are constants, and $x$ is an
auxiliary
independent variable (substituting any function $x=x(t)$ gives a solution of %
\eqref{evol4f} in terms of $t$).

Taking $f=f_1f_2f_3$ and all $f_i$ different, we obtain explicit expression of  the triaxial  Bianchi VIII
solutions indicated in \cite{BGPP}.

A particular solution is obtained by letting $a_1=a_3=0, a_2=\frac{a}{16}$ which gives
\begin{equation*}
f_1=f_3=\frac12(a-t^4)^{\frac14}, \quad f_2=\frac{t^2}2(a-t^4)^{-\frac14}, \quad f=t(a-t^4)^{-\frac14}, \quad
-a<t^4<a.
\end{equation*}
The resulting hyper K\"ahler metric is given by
\begin{equation*}
g=\frac12(a-t^4)^{\frac14}\Big((e^1)^2+\frac{t^2}{(a-t^4){^\frac12}}
(e^2)^2+(e^3)^2+\frac{2t}{(a-t^4)^{\frac12}}dt^2\Big),
\end{equation*}
where  the forms $e^i$ are given by \eqref{su11c}.

\medskip

\noindent
{\bf The Heisenberg group $H^3$, Bianchi type II, Gibbons-Hawking
class.}
Consider the two-step nilpotent Heisenberg group $H^3$ defined by the structure equations
\begin{equation}  \label{heis3}
\begin{aligned} de^1=de^2=0, \qquad de^3=-e^{12};\quad
e^1=dx, \quad e^2=dy, \quad e^3=dz-\frac12xdy+\frac12ydx. \end{aligned}
\end{equation}
The necessary conditions \eqref{4inst1} are satisfied. We evolve the structure according to \eqref{ev}. The
structure equations \eqref{heis3} reduce the evolution equations \eqref{evol4} to the following system of ODEs
\begin{equation}  \label{evol4h}
\frac{\partial}{\partial t}(f_1f_2)=ff_3,\quad \frac{\partial}{\partial t}%
(f_1f_3)=0, \quad \frac{\partial}{\partial t}(f_2f_3)=0.
\end{equation}
Working as in the previous example, i.e., using the same substitutions we see that the function $x_i$ satisfy
the system
\begin{equation*}
\frac {dx_3}{dr} = 2(x_1x_2x_3)^{1/4},\qquad \frac {dx_1}{dr} = \frac {dx_2}{%
dr}=0.
\end{equation*}
The general solution of thus system is
\begin{equation}  \label{e:H}
x_1=a, \qquad x_2=b, \qquad x_3=\left (\frac 32 (ab)^{1/4}\, r +c \right) ^{4/3},
\end{equation}
where $a, b$ and $c$ are constants. Therefore, using again  $f_i=\left ( \frac {x_jx_k}{x_i}\right )^{1/4} $,
the general solution of \eqref{evol4h} is
\begin{equation}  \label{e:H sols for xi}
\begin{aligned} f_1= \left (\frac ba \right )^{1/4}\left (\frac 32
(ab)^{1/4}\, r +c \right) ^{1/3}, \quad f_2= \left (\frac ab
\right )^{1/4}\left (\frac 32 (ab)^{1/4}\, r +c \right)
^{1/3},\quad f_3= \frac {\left (ab \right ) ^{1/4}}{\left (\frac
32 (ab)^{1/4}\, r +c \right ) ^{1/3}}. 
\end{aligned}
\end{equation}
A particular solution is obtained by taking $c=0$ and $a=b=1$, which gives
\begin{equation*}
f_1=f_2=\lambda r^{1/3}, \qquad f_3=f_1 ^{-1},
\end{equation*}
with $\lambda=\left (\frac 32 \right )^{1/3}$. The substitution $%
t=\lambda^2 r^{2/3}$ gives $f_1=f_2=f=t^{\frac12}, \quad f_3=t^{-\frac12}$. This is the hyper K\"ahler metric,
first written in \cite{Lor,LP1},
\begin{equation*}
g=t\Big[dt^2+dx^2+dy^2\Big]+\frac1t\Big[dz-\frac12xdy+\frac12ydx\Big]^2
\end{equation*}
belonging to the Gibbons-Hawking class \cite{GH} with an $S^1$-action and known also as Heisenberg metric
\cite{GR} (see also \cite{Bar,NP,SD,Cv04,VY}).

\medskip

\noindent {\bf Rigid motions of euclidean 2-plane-Bianchi $VII_0$ metrics.}
We consider the group $E_2$ of rigid motions of Euclidean 2-plane defined by the structure equations
\begin{equation}  \label{euc}
\begin{aligned} de^1=0, \quad de^2=e^{13},\quad de^3=-e^{12};\quad e^1=d\phi,
\quad e^2=\sin\phi dx-\cos\phi dy, \quad e^3 =\cos\phi dx+\sin\phi dy.
\end{aligned}
\end{equation}
Clearly \eqref{4inst1} are satisfied. We evolve the structure as in %
\eqref{ev}. Using the structure equations \eqref{euc} we reduce the evolution equations \eqref{evol4} to the
following system of ODE
\begin{equation}  \label{evol4feuc}
\frac{\partial}{\partial t}(f_1f_2)=ff_3,\quad \frac{\partial}{\partial t}%
(f_1f_3)=ff_2, \quad \frac{\partial}{\partial t}(f_2f_3)=0.
\end{equation}
With the substitution $x_{i}=(f_{j}f_{k})^{2}$, the above system becomes
\begin{equation*}
\frac{dx_{1}}{dr}=0,\qquad \frac{dx_{2}}{dr}=\frac{dx_{3}}{dr}=2(x_{1}x_{2}x_{3})^{1/4},
\end{equation*}%
in terms of the parameter $dr=fdt$. Hence, there is a function $x(r)$ and three constants $a_1, \ a_2, \ a_3$,
such that, $x(r)=x_{2}+a_{2}=x_{3}+a_{3}$, $x_1=a_1$. The equation for $x(r)$ is
\begin{equation}
\frac{d{x}}{dr}=2\left( a_{1}(x-a_{2})(x-a_{3})\right) ^{1/4},\ \text{%
i.e.,\ }dr=\frac{1}{2}\frac{1}{\left( a_{1}(x-a_{2})(x-a_{3})\right) ^{1/4}}dx.  \label{e:BIX12}
\end{equation}%
If we let $g(x)=\frac{1}{2}\left( a_{1}(x-a_{2})(x-a_{3})\right) ^{-1/4}$, and take into account
$x_{i}=(f_{j}f_{k})^{2}$,
we see from \eqref{evol4feuc} that the functions $f_{i}(x)$ satisfy
\begin{equation*}
\frac{d}{dx}\left( (x-a_{i})^{1/2}\right) =g(x)f_{i}, \quad i=2, \ 3.
\end{equation*}%
Solving for $f_{i}$ we show that the general solution of \eqref{evol4f} is
\begin{equation}\label{e:Emot}
\begin{aligned}
f_{1}(x)=\frac{(x-a_{2})^{1/4}(x-a_{3})^{1/4}}{a_{1}^{1/4}}%
,\qquad
f_{2}(x)=\frac{a_{1}^{1/4}(x-a_{3})^{1/4}}{(x-a_{2})^{1/4}},\qquad
f_{3}(x)=\frac{a_{1}^{1/4}(x-a_{2})^{1/4}}{(x-a_{3})^{1/4}},\\
\qquad f(t)= g\left ( x (t)\right )\,x'(t),\qquad
g(x)=\frac{1}{2}\left( (a_{1}(x-a_{2})(x-a_{3})\right) ^{-1/4},
\end{aligned}
\end{equation}
where $a_{1},\ a_{2}$ and $a_{3}$ are constants, and $x$ is an
auxiliary
independent variable (substituting any function $x=x(t)$ gives a solution of %
\eqref{evol4feuc} in terms of $t$).

\medskip

\noindent {\bf Vacuum solutions of Bianchi type $VII_0$.} When $f_2=f_3^{-1}, \quad f_1=f$ we have
\begin{equation*}
\frac{\partial}{\partial t}(ff_3^{-1})=ff_3,\qquad \frac{\partial}{\partial t%
}(ff_3)=ff_3^{-1},
\end{equation*}
with solution of the form $ff_3+ff_3^{-1}=Ae^t,\qquad ff_3^{-1}-ff_3=Be^{-t}$%
. Hence,
\begin{equation*}
f=f_1=\frac12(Ae^t+Be^{-t})^{\frac12}(Ae^t-Be^{-t})^{\frac12},\quad
f_3=f_2^{-1}=(Ae^t+Be^{-t})^{-\frac12}(Ae^t-Be^{-t})^{\frac12},
\end{equation*}
and the hyper K\"ahler metric is
\begin{equation}\label{hypeuc}
g=\frac14(A^2e^{2t}-B^2e^{-2t})\Big(dt^2+d\phi^2+%
\frac4{(Ae^t-Be^{-t})^2}(e^2)^2 +\frac4{(Ae^t+Be^{-t})^2}(e^3)^2\Big),
\end{equation}
where $e^2,e^3$ are given by \eqref{euc}.

In particular, setting $A=B$ in \eqref{hypeuc} we obtain
\begin{equation*}
g=\frac{A^2}2\sinh{2t}\Big(dt^2+d\phi^2\Big)+\coth t(e^2)^2 +\tanh
t(e^3)^2, \end{equation*} which is the vacuum solutions of Bianchi
type $VII_0$ \cite{Lor,LP1} with group of isometries $E_2$
\cite{GR}, (see also \cite{VY}).

\medskip

\noindent {\bf Rigid motions of Lorentzian 2-plane-Bianchi $VI_0$ metrics.}
Now we consider the group of rigid motions $E(1,1)$ of Lorentzian 2-plane defined by the structure equations
and coordinates as follows
\begin{equation}  \label{lor}
\begin{aligned} de^1=0, de^2=e^{13}, de^3=e^{12}; \quad 
e^1=d\phi, \quad e^2=\sinh\phi \, dx+\cosh\phi\, dy, \quad e^3 =\cosh\phi\, dx+\sinh\phi\, dy. \end{aligned}
\end{equation}
We evolve the structure as in \eqref{ev}. Using the structure equations %
\eqref{lor}, the evolution equations \eqref{evol4} turn into the next system of ODEs
\begin{equation}  \label{evol4flor}
\frac{\partial}{\partial t}(f_1f_2)=-ff_3,\quad \frac{\partial}{\partial t}%
(f_1f_3)=ff_2, \quad \frac{\partial}{\partial t}(f_2f_3)=0.
\end{equation}
The general solution of (\ref{evol4flor}) is
\begin{equation}\label{e:evol4flor}
\begin{aligned}
f_{1}(x)=\frac{(x-a_{2})^{1/4}(a_{3}-x)^{1/4}}{a_{1}^{1/4}}%
,\qquad
f_{2}(x)=\frac{a_{1}^{1/4}(a_{3}-x)^{1/4}}{(x-a_{2})^{1/4}},\qquad
f_{3}(x)=\frac{a_{1}^{1/4}(x-a_{2})^{1/4}}{(a_{3}-x)^{1/4}},\\
\qquad f(t)= g\left ( x (t)\right )\,x'(t),\quad
g(x)=\frac{1}{2}\left( (a_{1}(x-a_{2})(a_{3}-x)\right) ^{-1/4},
\end{aligned}
\end{equation}
where $a_{1},\ a_{2}$ and $a_{3}$ are constants, and $x$ is an
auxiliary
independent variable (substituting any function $x=x(t)$ gives a solution of %
\eqref{evol4feuc} in terms of $t$).

When $f_2=f_3^{-1}, \quad f_1=f$ we have $\frac{\partial}{\partial t}%
(ff_3^{-1})=-ff_3,\quad \frac{\partial}{\partial t}(ff_3)=ff_3^{-1} $ with solution of the form
\begin{equation*}
f=f_1=\frac12(a\cos t+b\sin t)^{\frac12}(a\cos t-b\sin t)^{\frac12},\quad f_3=f_2^{-1}=(a\cos t+b\sin
t)^{\frac12}(a\sin t-b\cos t)^{-\frac12},
\end{equation*}
and the hyper K\"ahler metric is given by
\begin{multline}  \label{hyplor}
g=\frac14(a^2\sin^2 t-b^2\cos^2 t)\Big(dt^2+d\phi^2+\frac4{(a\sin t+b\cos t)^2}(e^2)^2 +\frac4{(a\sin t-b\cos
t)^2}(e^3)^2\Big),
\end{multline}
where $e^2,e^3$ are given by \eqref{lor}. Introducing $t_0$ and $r_0$ by letting $r_0=\sqrt{a^2+b^2}$, $\cos
t_0=a/\sqrt{a^2+b^2}$ and $\sin t_0=b/\sqrt{a^2+b^2}$ the above metric can be put in the form
\begin{multline}  \label{e:hyplor}
g=\frac14(r_0^2\sin (t+t_0)\sin(t-t_0))\Big(dt^2+d\phi^2+\frac4{r_0^2\sin^2(t+t_0)}(e^2)^2 +
\frac4{r_0^2\sin^2(t-t_0)}(e^3)^2\Big).
\end{multline}

\medskip

\noindent {\bf Bianchi type $VI_0$} In particular, setting $a=b$ in \eqref{e:hyplor} we obtain
$r_0^2=2a^2, \quad \sin t_0=\cos t_0=\frac{\sqrt 2}{2}$. Taking $\tau=t+\frac{\pi}{4}$, the metric
\eqref{e:hyplor} takes the form
\begin{equation*}
g=\frac{a^2}4\sin{2\tau}\Big(d\tau^2+d\phi^2\Big)+\cot \tau(e^2)^2 +\tan \tau(e^3)^2,
\end{equation*}
which is the vacuum solutions of Bianchi type $VI_0$ \cite{Lor,LP1} with group of isometries $ E_(1,1)$
\cite{GR}, (see also \cite{VY}).

\section{Hyper symplectic (hyper para K\"ahler) metrics in dimension
4}

In this section, following the method of the preceding section, we present explicit hyper symplectic (hyper
para K\"ahler) metrics in dimension four of signature (2,2).
 The construction gives a kind of duality between hyper K\"ahler instantons and hyper para
K\"ahler structures.

We recall that an almost hyper paracomplex structure on a $4n$ dimensional space is a triple $(J,P_1,P_2)$
satisfying the paraquaternionic identities
$$J^2=-P_1^2=-P_2^2=-1, \quad JP_1=-P_1J=P_2.$$
A compatible metric $g$ satisfies
$$g(J.,J.)=-g(P_1.,P_1.)=-g(P_2.,P_2.)=g(.,.)$$
and is necessarily of neutral signature (2n,2n). The fundamental 2-forms are defined by
$$\Omega_1=g(.,J.),\quad \omega_2=g(.,.P_1), \quad \Omega_3=g(.,P_2).$$
When these forms are closed the structure is said to be hypersymplectic \cite{Hit1}. This implies (adapting the
computations of Atiyah-Hitchin \cite{AH} for hyper K\"ahler manifolds) that the structures are integrable and
parallel with respect to the Levi-Civita connection \cite{Hit1,DJS}.  Sometimes a hyper symplectic structure is
called also neutral hyper K\"ahler \cite{Kam}, hyper para K\"ahler \cite{IZ}. In dimension 4 an almost hyper
paracomplex structure is locally equivalent to an oriented neutral conformal structure, or an $Sp(1,\mathbb R)$
structure, and the integrability implies the anti-self-duality of the corresponding neutral conformal structure
\cite{Kam,IZ}. In particular, a hyper symplectic structure in dimension four underlines an anti-self-dual
Ricci-flat neutral metric. For this reason such structures  have been used in string theory
\cite{OV,hul,JR,Bar,Hull,CHO} and integrable systems \cite{D12,BM,DW}.

Let $G_3$ be a three dimensional Lie group with Lie algebra $g_3$ and $%
e^1,e^2,e^3$ be a basis of left invariant 1-forms. We consider the $Sp(1,\mathbb R)$ structure on
$g_3\times\mathbb{R}^+$ defined by the following 2-forms
\begin{equation}  \label{p4inst}
\begin{aligned} \Omega_1(t)=-e^1(t)\wedge e^2(t)+e^3(t)\wedge
f(t)dt,\\ \Omega_2(t)=e^1(t)\wedge e^3(t)-e^2(t)\wedge
f(t)dt,\\
\Omega_3(t)=e^2(t)\wedge e^3(t)+e^1(t)\wedge f(t)dt,
\end{aligned}
\end{equation}
where $f(t)$ is a function of $t$ and $e^i(t)$ depend on $t$.

With the help of Hitchin's theorem \cite{Hit1}, it is straightforward to prove the next
\begin{prop}
The $Sp(1,\mathbb R)$ structure $(\Omega_1,\Omega_2,\Omega_3)$ is hyper para K\"ahler if and only if
\begin{equation}\label{p4inst1} de^{12}=de^{13}=de^{23}=0,
\end{equation}
and the following evolution equations hold
\begin{equation}\label{pevol4}
\frac{\partial}{\partial t}e^{12}(t)=f(t)de^3(t),\quad \frac{\partial}{\partial t}e^{13}(t)=f(t)de^2(t), \quad
\frac{\partial}{\partial t}e^{23}(t)=-f(t)de^1(t).
\end{equation}
The hyper para K\"ahler metric is given by
\begin{equation}  \label{phypkel4}
g=(e^1)^2+(e^2)^2-(e^3)^2-f^2(t)dt^2.
\end{equation}
\end{prop}

\medskip

\noindent {\bf The group $SU(2)$.}
Let $G_3=SU(2)=S^3$ be described by the structure equations \eqref{su2}. Clearly \eqref{p4inst1} are satisfied.
We evolve the $SU(2)$ structure according to \eqref{ev}.

Using the structure equations \eqref{su2}, we reduce the evolution equations %
\eqref{pevol4} to the following system of ODEs
\begin{equation}  \label{pevol4f}
\frac{d}{d t}(f_1f_2)=-ff_3,\quad \frac d{dt}(f_1f_3)=ff_2, \quad \frac d{dt}(f_2f_3)=ff_1,
\end{equation}
which is equivalent to the system \eqref{evol4f11} after interchanging $f_2$ with $f_3$. The general solution
is given by \eqref{e:su2g}.

Taking $f=f_1f_2f_3$ in \eqref{e:su2g} and all $f_i$ different we obtain explicit expression of a triaxial
neutral hyper para K\"ahler metric
$$g=f_1^2(e_1)^2+f_3^2(e_2)^2-f_2^2(e_3)^2-f^2dt^2,
$$ where the forms $e^i$ are given by \eqref{eulersu2}.

A particular solution is obtained by letting $a_1=a_3=0, c_2=\frac{a}{16}$ in\eqref{e:su2g} which gives
\begin{equation*}
f_1=f_3=\frac12(a-r^4)^{\frac14}, \quad f_2=\frac{r^2}2(a-r^4)^{-\frac14}, \quad f=r(a-r^4)^{-\frac14}, \quad
-a<t^4<a.
\end{equation*}
The resulting neutral hyper para K\"ahler metric is
\begin{equation*}
g=\frac12(a-r^4)^{\frac14}\Big(d\theta^2+\sin^2\theta d\phi^2\Big)-\frac{r^2}{2(a-r^4){^\frac14}}
\Big(d\psi+\cos\theta d\phi\Big)^2-\frac{r}{(a-r^4)^{\frac14}}dr^2.
\end{equation*}

\medskip

\noindent {\bf The group $SU(1,1)$.}
Let $G_3=SU(1,1)$ be defined by the structure equations
\begin{equation}  \label{psu11}
de^1=-e^{23}, \quad de^2=-e^{31}, \quad de^3=e^{12}.
\end{equation}
In terms of local coordinates the left invariant forms $e^i$ are given by
\begin{equation}  \label{psu11c}
\begin{aligned}  e^1&=&d\psi
-\cos\theta d\phi,\quad e^2&=&\sinh\psi
d\theta+\cosh\psi\sin\theta d\phi, \quad e^3&=&\cosh\psi
d\theta+\sinh\psi\sin\theta d\phi. \end{aligned}
\end{equation}
Clearly \eqref{p4inst1} are satisfied. We consider the $SU(1,1)$ structure as in \eqref{ev}. Using the
structure equations \eqref{psu11c},  the evolution equations \eqref{pevol4} reduce to the already solved system
\eqref{evol4f} with a general solution of the form \eqref{e:su2}.

A particular solution to \eqref{evol4f} is given by \eqref{ehfff}, which results in a neutral hyper para
K\"ahler metric in Eguchi-Hanson form given by
\begin{multline*}
g=\frac{t^2}4\Big[\Big(d\psi -\cos\theta d\phi\Big)^2+\Big(\sinh\psi d\theta+\cosh\psi\sin\theta
d\phi\Big)^2\Big]\\-\frac{t^2}4\Big(1-\frac{a}{t^4}\Big)\Big(\cosh\psi d\theta+\sinh\psi\sin\theta
d\phi\Big)^2- \Big(1-\frac{a}{t^4}\Big)^{-1}(dt)^2.
\end{multline*}

Setting $f=-\frac{f_2}t$ one obtains another neutral hyper para K\"ahler. Triaxial neutral hyper para K\"ahler
metric can be obtained with the help of \eqref{triax}.

\medskip

\noindent {\bf The Heisenberg group $H^3$.}
Consider the two-step nilpotent Heisenberg group $H^3$ defined by the structure equations \eqref{heis3}. The
structure equations \eqref{heis3} reduce the evolution equations \eqref{pevol4} to the already solved system
\eqref{evol4h} with a general solution \eqref{e:H sols for xi}.

A particular solution is  $f_1=f_2=f=t^{\frac12}, \quad f_3=-t^{-\frac12}$. This is the neutral hyper para
K\"ahler metric
\begin{equation*}
g=t\Big[-dt^2+dx^2+dy^2\Big]-\frac1t\Big[dz-\frac12xdy+\frac12ydx\Big]^2.
\end{equation*}

\medskip

\noindent {\bf Rigid motions of the Euclidean 2-plane.}
We consider the group $E_2$ of rigid motions of Euclidean 2-plane defined by the structure equations
\eqref{euc}. Clearly \eqref{4inst1} are satisfied. We evolve the structure as in \eqref{ev}. Using the
structure equations \eqref{euc}, the evolution equations \eqref{pevol4} take the form of the already solved
system of ODEs \eqref{evol4flor} with a general solution \eqref{e:evol4flor}.

When $f_2=f_3^{-1}, \quad f_1=f$ we have
\begin{equation*}
f=f_1=\frac12(a\cos t+b\sin t)^{\frac12}(a\cos t-b\sin t)^{\frac12},\quad f_3=f_2^{-1}=(a\cos t+b\sin
t)^{\frac12}(a\sin t-b\cos t)^{-\frac12}.
\end{equation*}
Introducing $t_0$ and $r_0$ by letting $r_0=\sqrt{a^2+b^2}$, $\cos t_0=a/\sqrt{a^2+b^2}$ and $\sin
t_0=b/\sqrt{a^2+b^2}$, the resulting neutral hyper para K\"ahler metric can be put in the form
\begin{multline}  \label{e:phyplor}
g=\frac14(r_0^2\sin (t+t_0)\sin(t-t_0))\Big(-dt^2+d\phi^2+\frac4{r_0^2\sin^2(t+t_0)}(e^2)^2 -
\frac4{r_0^2\sin^2(t-t_0)}(e^3)^2\Big),
\end{multline}
where $e^2,e^3$ are given by \eqref{euc}.

In particular, setting $a=b$ in \eqref{e:hyplor} we obtain $r_0^2=2a^2, \quad \sin t_0=\cos t_0=\frac{\sqrt
2}{2}$. Taking $\tau=t+\frac{\pi}{4}$, the metric \eqref{e:phyplor}
can be written as
\begin{equation*}
g=\frac{a^2}4\sin{2\tau}\Big(-d\tau^2+d\phi^2\Big)+\cot \tau\Big(\sin\phi \, dx-\cos\phi\, dy\Big)^2 -\tan
\tau\Big(\cos\phi\, dx+\sin\phi\, dy\Big)^2.
\end{equation*}

\medskip

\noindent {\bf Rigid motions of Lorentzian 2-plane-Bianchi $VI_0$ metrics.}
Now we consider the group of rigid motions $E(1,1)$ of Lorentzian 2-plane defined by the structure equations
\eqref{lor}.
We evolve the structure as in \eqref{ev}. Using the structure equations %
\eqref{lor}, the evolution equations \eqref{pevol4} turn into the solved system of ODEs \eqref{evol4feuc} with
the general solution given by \eqref{e:Emot}.

When $f_2=f_3^{-1}, \quad f_1=f$ we have
\begin{equation*}
f=f_1=\frac12(Ae^t+Be^{-t})^{\frac12}(Ae^t-Be^{-t})^{\frac12},\quad
f_3=f_2^{-1}=(Ae^t+Be^{-t})^{-\frac12}(Ae^t-Be^{-t})^{\frac12},
\end{equation*}
and the neutral hyper para K\"ahler metric is
\begin{multline}  \label{phypeuc}
g=\frac14(A^2e^{2t}-B^2e^{-2t})\Big(-dt^2+d\phi^2+%
\frac4{(Ae^t-Be^{-t})^2}(e^2)^2 -\frac4{(Ae^t+Be^{-t})^2}(e^3)^2\Big),
\end{multline}
where $e^2,e^3$ are given by \eqref{lor}.

In particular, setting $A=B$ in \eqref{phypeuc} we obtain
\begin{equation*}
g=\frac{A^2}2\sinh{2t}\Big(-dt^2+d\phi^2\Big)+\coth t\Big( \sinh\phi \, dx+\cosh\phi\, dy\Big)^2 -\tanh t\Big(
\cosh\phi \, dx+\sinh\phi\, dy\Big)^2.
\end{equation*}

\section{Hyper K\"ahler structures in dimension eight}\label{h8kel}

In this section we apply our method from Section~\ref{sph}.

Let $G_7$ be the seven dimensional solvable non-nilpotent Lie group defined by the following structure
equations
\begin{equation}  \label{ex78}
\begin{aligned}
& de^1=e^{17}+e^{27},\quad de^2=-e^{17}-e^{27},\quad
de^3=-e^{15}+e^{16}-e^{25}+e^{26},\quad
 de^4=-e^{16}-e^{15}-e^{25}-e^{26},\\ &
 de^5=e^{13}+e^{14}+e^{23}+e^{24},\qquad
 de^6=-e^{13}+e^{14}-e^{23}+e^{24},\qquad de^7=2 e^{12}.
\end{aligned}
\end{equation}
This is a solvable non nilpotent Lie algebra because $[g,g] = g_1$ is
generated by ${e_1-e_2, e_3, e_4, e_5, e_6, e_7}$, $[g,g_1]=g_1$ and $%
[g_1,g_1]=0$. The $Sp(2)$-hypo structure is determined by the equalities
\begin{equation*}
\begin{aligned}
 d(e^{12}+e^{34}+e^{56})=0,\qquad
     d(e^{13}-e^{24}+e^{57})=0,\qquad
    d(e^{14}+e^{23}+e^{67})=0.
    \end{aligned}
\end{equation*}

We consider the $Sp(2)$ structure on $g_3\times\mathbb{R}^+$ defined by the following 2-forms
\begin{equation}  \label{4inst8}
\begin{aligned} F_1(t)=e^1(t)\wedge e^2(t)+e^3(t)\wedge e^4(t)+ e^5(t)\wedge
e^6(t)+e^7(t)\wedge f(t)dt,\\ F_2(t)=e^1(t)\wedge e^3(t)-e^2(t)\wedge
e^4(t)+ e^5(t)\wedge e^7(t)-e^6(t)\wedge f(t)dt,\\
F_3(t)=e^1(t)\wedge e^4(t)+e^2(t)\wedge e^3(t)+ e^6(t)\wedge e^7(t)+e^5(t)\wedge f(t)dt.
\end{aligned}
\end{equation}
where $f(t)$ is a function of $t$ and $e^i(t)$ depend on $t$. A direct calculation shows that for the evolution
\begin{equation}\label{lie8}
e^1(t)=-te^1-(t+1)e^2, \qquad e^2(t)=-(t+1)e^1-te^2, \qquad e^a(t)=e^a, \ a=3,\dots, 7,
\end{equation}
the corresponding forms $F_1(t), \ F_2(t), \ F_3(t)$ are closed.

We consider the basis
\begin{equation}\label{newbas}
\epsilon^1=\sqrt{2}(e^1+e^2), \quad \epsilon^2=e^2,\quad
\epsilon^3=e^3+e^4,\quad \epsilon^4=e^3-e^4,\quad
\epsilon^5=\sqrt{2}e^5, \quad \epsilon^6=\sqrt{2}e^6,\quad
\epsilon^7=\frac1{\sqrt{2}}e^7.
\end{equation}
In this basis the structure equations \eqref{ex78} take the form
\begin{equation}  \label{ex78n}
\begin{aligned} d\epsilon^1=0,\quad d\epsilon^2=-\epsilon^{17},\quad
d\epsilon^3=-\epsilon^{15},\quad d\epsilon^4=\epsilon^{16},\quad
d\epsilon^5=\epsilon^{13},\quad d\epsilon^6=-\epsilon^{14},\quad
d\epsilon^7=\epsilon^{12}. \end{aligned}
\end{equation}
Considering the triples $(\epsilon^1,\epsilon^2,\epsilon^7)$, $(\epsilon^1,\epsilon^3,\epsilon^5)$,
$(\epsilon^1,\epsilon^4,\epsilon^6)$, we obtain
\begin{equation}\label{invf}
\begin{aligned}
& \epsilon^1=\ dx^1, \quad \epsilon^2=\cos x^1\,dx^2-\sin {x^1}\,dx^7,\quad \epsilon^7=
(\sin {x^1}\,dx^2+\cos {x^1}\,dx^7),\\
& \epsilon^3=-(\sin {x^1}\,dx^5+\cos {x^1}\,dx^3),\quad \epsilon^5=
\cos {x^1}\,dx^5-\sin {x^1}\,dx^3,\\
& \epsilon^4=(\sin {x^1}\,dx^6+\cos {x^1}\,dx^4),\quad \epsilon^6= \cos {x^1}\,dx^6 - \sin {x^1}\,dx^4.
\end{aligned}
\end{equation}
For the  hyper K\"ahler metric on $G_7\times\mathbb R$  given by
$g=\sum_{r=1}^7e^r(t)^2+dt^2 $
the equations \eqref{newbas} and \eqref{invf} yield 
\begin{equation*}\label{hkelfin}
g=(t^2+t+1/2)(dx^1)^2+2(dx^2)^2+2(dx^7)^2-\sqrt 2\cos x^1\,
dx^1dx^2+\sqrt 2\sin x^1\, dx^1dx^7+\sum_{s=3}^6(dx^s)^2+dt^2.
\end{equation*}
When $t=-1/2$ the metric degenerates ($e_1-e_2$ is of zero length). The above metric is of constant zero
curvature, but it is not complete.
 The $8$-dimensional manifold becomes a product of the Euclidean $\mathbb{R}^4$ with a four dimensional
 manifold $M$ of vanishing curvature.

One can consider also the following  $Sp(2)$ structure on $G_7\times\mathbb{R}^+$
\begin{equation}  \label{4inst8n}
\begin{aligned} F^1(t)=\epsilon^1(t)\wedge \epsilon^2(t)+\epsilon^3(t)\wedge \epsilon^4(t)- \epsilon^5(t)\wedge
\epsilon^6(t)+\epsilon^7(t)\wedge h(t)dt,\\ F^2(t)=\epsilon^1(t)\wedge \epsilon^3(t)-\epsilon^2(t)\wedge
\epsilon^4(t)- \epsilon^6(t)\wedge \epsilon^7(t)+\epsilon^5(t)\wedge h(t)dt,\\
F^3(t)=\epsilon^1(t)\wedge \epsilon^4(t)+\epsilon^2(t)\wedge \epsilon^3(t)- \epsilon^5(t)\wedge
\epsilon^7(t)-\epsilon^6(t)\wedge h(t)dt,
\end{aligned}
\end{equation}
where $h(t)$ is a function of $t$ and $\epsilon^i(t)$ depend on $t$. A direct calculation shows that for the
evolution
\begin{equation}\label{lie8n}
\epsilon^1(t)=h_1(t)\epsilon^1, \qquad \epsilon^a(t)=\epsilon^a, \ a=2,\dots, 7, \qquad h_1'=-h
\end{equation}
the corresponding forms $F^1(t), \ F^2(t), \ F^3(t)$ are closed. The corresponding hyper K\"ahler metric
$g=\sum_{r=1}^7\epsilon^r(t)^2+dt^2 $ is flat having the expression ($u=h_1(t)$)
$$
g=u^2(dx^1)^2+(du)^2+(dx^2)^2+(dx^3)^2+(dx^4)^2+(dx^5)^2+(dx^6)^2+(dx^7)^2.
$$
\newpage
\section{Appendix 1. Explicit quaternionic K\"ahler metric}
Substituting in \eqref{qknew8} the equations \eqref{loccoord} we obtain  the following expression for the
metric coefficients of the quaternionic K\"ahler metric \eqref{qknew8} in coordinates\\
$\{v^1=t,v^2=x,v^3=y,v^4=z,v^5=x_5,v^6=x_6,v^7=x_7,v^8=u\}$:

\[
\begin{array}{l}
g_{11}=\frac{1}{4} u \left(a u \left(x_5^2+x_6^2+x_7^2\right)+4\right),\qquad g_{12}=-\frac{1}{2} u (a u-1) x_5,\\[7pt]
g_{13}=\frac{1}{2} u (a u-1) \left(x_6\cos x -x_7\sin x \right),\\[7pt]
 g_{14}=-\frac{1}{2} u (a u-1) \left(x_5\cos y -\sin y \left(x_6\sin x+x_7\cos x\right)\right),
 \quad g_{15}=-\frac{1}{4} a u^2 x_5\\[7pt]
 g_{16}=-\frac{1}{4} a u^2 x_6,\quad g_{17}=-\frac{1}{4} a u^2 x_7,\\[7pt]
g_{22}=\frac{1}{4} u \left(a u \left(x_6^2+x_7^2+4\right)-4\right),\quad  g_{23}=\frac{1}{4} a u^2 x_5
\left(x_6\cos x-x_7\sin x\right),\\[7pt]
 g_{24}=\frac{1}{4} u \left(a u x_5\sin y \left(x_6\sin x+x_7\cos x\right)+\cos y \left(a u
   \left(x_6^2+x_7^2+4\right)-4\right)\right),\\[7pt]
g_{25}=\frac{1}{2} u (a u-1),\qquad g_{26}=-\frac{1}{4} a u^2 x_7,\qquad g_{27}=\frac{1}{4} a u^2 x_6,\\[7pt]
g_{33}=\frac{1}{8} u\left(2 a ux_5^2+a ux_6^2+a ux_7^2+8 a u+2 a ux_6x_7\sin 2 x\right.\\[5pt]
    \qquad\ \left.-a u \cos 2x\left(x_6^2-x_7^2\right)-8\right),\\[7pt]
g_{34}=\frac{1}{8} a u^2 \left(x_6\cos x-x_7\sin x\right) \left(2x_5 \cos y-2\left(x_6\sin x+x_7\cos x\right)
\sin y \right),\\[7pt]
g_{35}=-\frac{1}{4} a u^2 \left(x_6\sin x+x_7\cos x\right),\qquad g_{36}=\frac{1}{4} u \bigl((2-2 a u)
\cos x+a u x_5\sin x\bigr),\\[7pt]
 g_{37}=\frac{1}{4} u \bigl(2 (a u-1) \sin x+a ux_5\cos x\bigr),\\[7pt]
 g_{44}=\frac{1}{4} u \left(\left(a u \left(x_6^2+x_7^2+4\right)-4\right) \cos ^2y+4 a u \sin ^2x \sin ^2y
 -4 \sin ^2x\sin ^2y\right.\\[5pt]
 \qquad\ \left.
 +a ux_5^2 \sin ^2x \sin ^2x +a ux_7^2\sin ^2x \sin ^2x +a u x_5x_6 \sin x \sin 2y \right.\\[5pt]
 \qquad\ \left.+\cos^2x \sin ^2y \left(a u \left(x_5^2+x_6^2+4\right)-4\right)+a u x_5x_7\cos x \sin 2y\right.\\[5pt]
 \qquad\ \left. -a ux_6x_7 \sin 2x \sin^2y\right),\\[7pt]
 g_{45}=\frac{1}{4} u \left(2 (a u-1) \cos y+a u \sin y \left(x_6\cos x-x_7\sin x\right)\right),\\[7pt]
 g_{46}=-\frac{1}{4} u \left(2 (a u-1) \sin x \sin y+a u \left(x_5\cos x \sin y +x_7\cos (y)\right)\right),\\[7pt]
 g_{47}=\frac{1}{4} u \left(a u \left(x_5\sin x \sin y+x_6\cos y\right)-2 (a u-1) \cos x \sin y\right),\\[7pt]
 g_{55}=g_{66}=g_{77}=\frac{a u^2}{4},\qquad g_{88}=\frac{1}{u (a u-1)}.
\end{array}
\]

\newpage

\section{Appendix 2. Explicit $Spin(7)$-holonomy metric}
Substituting in \eqref{newspin7metric} the equations \eqref{loccoord} we obtain  the following
expression for the metric coefficients of the $Spin(7)$ metric \eqref{newspin7metric} in coordinates\\
$\{v^1=t,v^2=x,v^3=y,v^4=z,v^5=x_5,v^6=x_6,v^7=x_7,v^8=u\}$:

\[
\begin{array}{l}
g_{11}=\frac{20 u^{5/3}+\left(9 u^{5/3}+4 a\right) \left(x_5^2+x_6^2+x_7^2\right)}{20 u^{2/3}},\qquad
g_{12}=-\frac{2 \left(u^{5/3}+a\right) x_5}{5 u^{2/3}},\\[7pt]
g_{13}=\frac{2 \left(u^{5/3}+a\right) \left(x_6\cos x-x_7\sin x\right)}{5 u^{2/3}},\quad
g_{14}=-\frac{2 \left(u^{5/3}+a\right) \left(x_5\cos y-\sin y\left(x_6\sin x+x_7\cos x\right)\right)}{5 u^{2/3}},\\[7pt]
g_{15}=-\frac{\left(9 u^{5/3}+4 a\right) x_5}{20 u^{2/3}},\qquad g_{16}=-\frac{\left(9 u^{5/3}+4 a\right) x_6}
{20 u^{2/3}},\qquad g_{17}=-\frac{\left(9 u^{5/3}+4 a\right) x_7}{20 u^{2/3}},\\[7pt]
g_{22}=\frac{16 \left(u^{5/3}+a\right)+\left(9 u^{5/3}+4 a\right) \left(x_6^2+x_7^2\right)}{20 u^{2/3}},\qquad
g_{23}=\frac{\left(9 u^{5/3}+4 a\right) x_5 \left(x_6\cos x-x_7\sin x\right)}{20 u^{2/3}},\\[7pt]
g_{24}=\frac{16 \left(u^{5/3}+a\right) \cos y+\left(9 u^{5/3}+4 a\right) \left(x_5\sin y \left(x_6\sin x+
x_7\cos x\right)+\left(x_6^2+x_7^2\right)\cos y \right)}{20 u^{2/3}},\\[7pt]
g_{25}=\frac{2 \left(u^{5/3}+a\right)}{5 u^{2/3}},\qquad g_{26}=-\frac{\left(9 u^{5/3}+4 a\right) x_7}
{20 u^{2/3}}, \qquad g_{27}=\frac{\left(9 u^{5/3}+4 a\right) x_6}{20 u^{2/3}},\\[7pt]
g_{33}=\frac{16 \left(u^{5/3}+a\right)+\left(9 u^{5/3}+4 a\right) \left(x_5^2+x_6^2+x_7^2-\left(x_6\cos x-x_7
\sin x\right)^2\right)}{20 u^{2/3}},\\[7pt]
g_{34}=-\frac{\left(9 u^{5/3}+4 a\right) \left(x_6\cos x-x_7\sin x\right) \left(-x_5\cos y+x_6\sin x \sin
y+x_7\cos
  x \sin (y\right)}{20 u^{2/3}},\\[7pt]
g_{35}=-\frac{\left(9 u^{5/3}+4 a\right) \left(x_6\sin x+x_7\cos x\right)}{20 u^{2/3}},\quad
g_{36}=\frac{\left(9 u^{5/3}+4 a\right)x_5\sin x-8 \left(u^{5/3}+a\right) \cos x}{20 u^{2/3}},\\[7pt]
g_{37}=\frac{8 \left(u^{5/3}+a\right) \sin x+\left(9 u^{5/3}+4 a\right)x_5 \cos x}{20 u^{2/3}},\\[7pt]
g_{44}=\frac{4 \left(u^{5/3}+a\right)}{5 u^{2/3}}-\frac{\left(9 u^{5/3}+4 a\right)x_5 \sin 2y
\left(x_6\sin x+x_7\cos x\right)}{20 u^{2/3}}\\[5pt]
\qquad\ -\frac{\left(9 u^{5/3}+4 a\right) \left(\left(x_6^2+x_7^2\right) \cos ^2y+ \left(x_5^2+
\left(x_6\cos x-x_7\sin x\right)^2\right)\sin ^2y\right)}{20 u^{2/3}},\\[7pt]
g_{45}=\frac{2 \left(u^{5/3}+a\right) \cos y}{5 u^{2/3}}+\frac{\left(9 u^{5/3}+4 a\right) \sin y
\left(x_6\cos x-x_7\sin x\right)}{20 u^{2/3}},\\[7pt]
g_{46}=-\frac{2 \left(u^{5/3}+a\right) \sin x \sin y}{5 u^{2/3}}-\frac{\left(9 u^{5/3}+4 a\right)
\left(x_5\cos x \sin y+x_7\cos y\right)}{20 u^{2/3}},\\[7pt]
g_{47}=\frac{\left(9 u^{5/3}+4 a\right) \left(x_5\sin x \sin y+x_6\cos y\right)}{20 u^{2/3}}-\frac{2
   \left(u^{5/3}+a\right) \cos x \sin y}{5 u^{2/3}},\\[7pt]
g_{55}=g_{66}=g_{77}=\frac{9 u^{5/3}+4 a}{20 u^{2/3}},\qquad g_{88}=\frac{5 u^{2/3}}{36 \left(u^{5/3}+a\right)}
\end{array}
\]

\subsection{Holonomy of the $Spin(7)$ metrics}

Let us consider the Lie group \eqref{ex11} and the metric
$$
g=u\left ((e^1)^2+(e^2)^2+(e^3)^2+(e^4)^2 \right )+ \frac{(a-u^{5/3})}{20u^{2/3}}\left (
(\eta_1)^2+(\eta_2)^2+(\eta_3)^2 \right )+ \frac{5u^{2/3}}{9(a-u^{5/3})}du^2.
$$
Since $\eta_1=e^5$, $\eta_2=e^6$ and $\eta_3=e^7$, the metric can be written as
$$
g=(\sqrt{u}\, e^1)^2+(\sqrt{u}\, e^2)^2+(\sqrt{u}\, e^3)^2+(\sqrt{u}\, e^4)^2 + (g(u)\, e^5)^2+(g(u)\,
e^6)^2+(g(u)\, e^7)^2 + \left( \frac{1}{6\, g(u)} du \right)^2,
$$
where the function $g(u)$ is given by $g(u)=\sqrt{\frac{(a-u^{5/3})}{20\, u^{2/3}}}$. From now on, we
shall work with the orthonormal basis
$$
\{\gamma^1=\sqrt{u}\, e^1, \gamma^2=\sqrt{u}\, e^2, \gamma^3=\sqrt{u}\, e^3, \gamma^4=\sqrt{u}\, e^4,
\gamma^5=g(u)\, e^5, \gamma^6=g(u)\, e^6, \gamma^7=g(u)\, e^7, \gamma^8=\frac{du}{6\, g(u)}\}.
$$
The curvature 2-forms $\Omega^i_j$ of the Levi-Civita connection with respect to the basis
$\{\gamma^1,\ldots,\gamma^8\}$ are:

\bigskip

$\Omega^1_2= -\frac{1\, u+12\, g(u)^2}{u^2}\gamma^{12} -\frac{6g(u)(2ug'(u)-g(u))}{u^2}\gamma^{58} -
\frac{1\, u+4\, g(u)^2}{2u^2}\gamma^{67} $

\bigskip

$\Omega^1_3= -\frac{1\, u+12\, g(u)^2}{u^2}\gamma^{13}+\frac{1\, u+4\, g(u)^2}{2u^2}\gamma^{57}
+\frac{6g(u)(2ug'(u)-g(u))}{u^2}\gamma^{68} $

\bigskip

$\Omega^1_4= -\frac{1\, u+12\, g(u)^2}{u^2}\gamma^{14} - \frac{1\, u+4\, g(u)^2}{2u^2}\gamma^{56}
+\frac{6g(u)(2ug'(u)-g(u))}{u^2}\gamma^{78} $

\bigskip

$\Omega^1_5= -\frac{g(u)(18ug'(u)-g(u))}{u^2} \gamma^{15} +\frac{3g(u)(2ug'(u)-g(u))}{u^2}\gamma^{28} +
\frac{1\, u+4\, g(u)^2}{4u^2}\gamma^{37} - \frac{1\, u+4\, g(u)^2}{4u^2}\gamma^{46} $

\bigskip

$\Omega^1_6= -\frac{g(u)(18ug'(u)-g(u))}{u^2} \gamma^{16} - \frac{1\, u+4\, g(u)^2}{4u^2}\gamma^{27}
+\frac{3g(u)(2ug'(u)-g(u))}{u^2}\gamma^{38} + \frac{1\, u+4\, g(u)^2}{4u^2}\gamma^{45} $

\bigskip

$\Omega^1_7= -\frac{g(u)(18ug'(u)-g(u))}{u^2} \gamma^{17} + \frac{1\, u+4\, g(u)^2}{4u^2}\gamma^{26} -
\frac{1\, u+4\, g(u)^2}{4u^2}\gamma^{35} +\frac{3g(u)(2ug'(u)-g(u))}{u^2}\gamma^{48} $

\bigskip

$\Omega^1_8= -\frac{9g(u)(2ug'(u)-g(u))}{u^2} \gamma^{18} -\frac{3g(u)(2ug'(u)-g(u))}{u^2}\gamma^{25} -
\frac{3g(u)(2ug'(u)-g(u))}{u^2}\gamma^{36} - \frac{3g(u)(2ug'(u)-g(u))}{u^2}\gamma^{47} $

\bigskip

$\Omega^2_3= -\frac{1\, u+12\, g(u)^2}{u^2} \gamma^{23} -\frac{1\, u+4\, g(u)^2}{2u^2} \gamma^{56} +
\frac{6g(u)(2ug'(u)-g(u))}{u^2}\gamma^{78} $

\bigskip

$\Omega^2_4= -\frac{1\, u+12\, g(u)^2}{u^2} \gamma^{24} -\frac{1\, u+4\, g(u)^2}{2u^2} \gamma^{57} -
\frac{6g(u)(2ug'(u)-g(u))}{u^2}\gamma^{68} $

\bigskip

$\Omega^2_5= - \frac{3g(u)(2ug'(u)-g(u))}{u^2}\gamma^{18} - \frac{g(u)(18ug'(u)-g(u))}{u^2}\gamma^{25}
-\frac{1\, u+4\, g(u)^2}{4u^2} \gamma^{36} -\frac{1\, u+4\, g(u)^2}{4u^2} \gamma^{47} $

\bigskip

$\Omega^2_6= \frac{1\, u+4\, g(u)^2}{4u^2} \gamma^{17} - \frac{g(u)(18ug'(u)-g(u))}{u^2}\gamma^{26}
+\frac{1\, u+4\, g(u)^2}{4u^2} \gamma^{35} - \frac{3g(u)(2ug'(u)-g(u))}{u^2}\gamma^{48} $

\bigskip

$\Omega^2_7= -\frac{1\, u+4\, g(u)^2}{4u^2} \gamma^{16} - \frac{g(u)(18ug'(u)-g(u))}{u^2}\gamma^{27} +
\frac{3g(u)(2ug'(u)-g(u))}{u^2}\gamma^{38} +\frac{1\, u+4\, g(u)^2}{4u^2} \gamma^{45} $

\bigskip

$\Omega^2_8= \frac{3g(u)(2ug'(u)-g(u))}{u^2} \gamma^{15} -\frac{9g(u)(2ug'(u)-g(u))}{u^2} \gamma^{28}
-\frac{3g(u)(2ug'(u)-g(u))}{u^2} \gamma^{37} + \frac{3g(u)(2ug'(u)-g(u))}{u^2} \gamma^{46} $

\bigskip

$\Omega^3_4= -\frac{1\, u+12\, g(u)^2}{u^2} \gamma^{34} +\frac{6g(u)(2ug'(u)-g(u))}{u^2} \gamma^{58}
-\frac{1\, u+4\, g(u)^2}{2u^2} \gamma^{67} $

\bigskip

$\Omega^3_5= -\frac{1\, u+4\, g(u)^2}{4u^2} \gamma^{17}+\frac{1\, u+4\, g(u)^2}{4u^2} \gamma^{26}
-\frac{g(u)(18ug'(u)-g(u))}{u^2} \gamma^{35} + \frac{3g(u)(2ug'(u)-g(u))}{u^2} \gamma^{48} $

\bigskip

$\Omega^3_6=  - \frac{3g(u)(2ug'(u)-g(u))}{u^2} \gamma^{18} - \frac{1\, u+4\, g(u)^2}{4u^2} \gamma^{25}
-\frac{g(u)(18ug'(u)-g(u))}{u^2} \gamma^{36} -\frac{1\, u+4\, g(u)^2}{4u^2} \gamma^{47} $

\bigskip

$\Omega^3_7= \frac{1\, u+4\, g(u)^2}{4u^2} \gamma^{15} - \frac{3g(u)(2ug'(u)-g(u))}{u^2} \gamma^{28}
-\frac{g(u)(18ug'(u)-g(u))}{u^2} \gamma^{37} +\frac{1\, u+4\, g(u)^2}{4u^2} \gamma^{46} $

\bigskip

$\Omega^3_8= \frac{3g(u)(2ug'(u)-g(u))}{u^2} \gamma^{16} +\frac{3g(u)(2ug'(u)-g(u))}{u^2} \gamma^{27}
-\frac{9g(u)(2ug'(u)-g(u))}{u^2} \gamma^{38} - \frac{3g(u)(2ug'(u)-g(u))}{u^2} \gamma^{45} $

\bigskip

$\Omega^4_5= \frac{1\, u+4\, g(u)^2}{4u^2} \gamma^{16}+\frac{1\, u+4\, g(u)^2}{4u^2} \gamma^{27}
-\frac{3g(u)(18ug'(u)-g(u))}{u^2} \gamma^{38} - \frac{g(u)(18ug'(u)-g(u))}{u^2} \gamma^{45} $

\bigskip

$\Omega^4_6= - \frac{1\, u+4\, g(u)^2}{4u^2} \gamma^{15} + \frac{3g(u)(2ug'(u)-g(u))}{u^2} \gamma^{28}
+\frac{1\, u+4\, g(u)^2}{4u^2} \gamma^{37} - \frac{g(u)(18ug'(u)-g(u))}{u^2} \gamma^{46} $

\bigskip

$\Omega^4_7= - \frac{3g(u)(2ug'(u)-g(u))}{u^2} \gamma^{18} - \frac{1\, u+4\, g(u)^2}{4u^2} \gamma^{25} -
\frac{1\, u+4\, g(u)^2}{4u^2} \gamma^{36} -\frac{g(u)(18ug'(u)-g(u))}{u^2} \gamma^{47} $

\bigskip

$\Omega^4_8= \frac{3g(u)(2ug'(u)-g(u))}{u^2} \gamma^{17} -\frac{3g(u)(2ug'(u)-g(u))}{u^2} \gamma^{26} +
\frac{3g(u)(2ug'(u)-g(u))}{u^2} \gamma^{35} -\frac{9g(u)(2ug'(u)-g(u))}{u^2} \gamma^{48} $

\bigskip

$\Omega^5_6= - \frac{1\, u+4\, g(u)^2}{2u^2} \gamma^{14} - \frac{1\, u+4\, g(u)^2}{2u^2} \gamma^{23} -
\frac{(24\, g(u) g'(u)+\lambda^2 +\mu^2)(24\, g(u) g'(u)-\lambda^2 -\mu^2)}{16\, g(u)^2} \gamma^{56} $

\bigskip

$\Omega^5_7= \frac{1\, u+4\, g(u)^2}{2u^2} \gamma^{13} - \frac{1\, u+4\, g(u)^2}{2u^2} \gamma^{24} -
\frac{(24\, g(u) g'(u)+\lambda^2 +\mu^2)(24\, g(u) g'(u)-\lambda^2 -\mu^2)}{16\, g(u)^2} \gamma^{57} $

\bigskip

$\Omega^5_8= \frac{6g(u)(2ug'(u)-g(u))}{u^2} \gamma^{12} + \frac{6g(u)(2ug'(u)-g(u))}{u^2} \gamma^{34} -
36(g'(u)^2+g(u)g''(u)) \gamma^{58} $

\bigskip

$\Omega^6_7= -\frac{1\, u+4\, g(u)^2}{2u^2} \gamma^{12} - \frac{1\, u+4\, g(u)^2}{2u^2} \gamma^{34} -
\frac{(24\, g(u) g'(u)+\lambda^2 +\mu^2)(24\, g(u) g'(u)-\lambda^2 -\mu^2)}{16\, g(u)^2} \gamma^{67} $

\bigskip

$\Omega^6_8= \frac{6g(u)(2ug'(u)-g(u))}{u^2} \gamma^{13} - \frac{6g(u)(2ug'(u)-g(u))}{u^2} \gamma^{24} -
36(g'(u)^2+g(u)g''(u)) \gamma^{68} $

\bigskip

$\Omega^7_8= \frac{6g(u)(2ug'(u)-g(u))}{u^2} \gamma^{14} + \frac{6g(u)(2ug'(u)-g(u))}{u^2} \gamma^{23} -
36(g'(u)^2+g(u)g''(u)) \gamma^{78}.$

\bigskip

First of all, using that $g(u)=\sqrt{\frac{1(a-u^{5/3})}{20\, u^{2/3}}}$, from these expressions one can
check directly that the metric is Ricci flat because
$$
{\rm Ric}(X_i,X_j)= \Omega^1_j(X_1,X_i)+\cdots+\Omega^8_j(X_8,X_i) =0,
$$
for any $i,j=1,\ldots,8$ and for any $a$, where $\{X_1,\ldots,X_8\}$ denotes the dual basis of
$\{\gamma^1,\ldots,\gamma^8\}$. Now, one can evaluate the coefficients above using that
$g(u)=\sqrt{\frac{1(a-u^{5/3})}{20\, u^{2/3}}}$. It turns out that all the coefficients above are nonzero
when $a\not=0$ and $\lambda^2+\mu^2\not=0$. It is clear that the first 9 curvature forms, i.e from $\Omega^1_2$
to $\Omega^2_4$, are independent. The form $\Omega^2_5$ is independent from the previous ones, except possibly
for $\Omega^1_8$. But if $\Omega^1_8$ and $\Omega^2_5$ were proportional then, from the coefficient in
$\gamma^{18}$, the factor of proportionality should be equal to 3 and this is not the case for the coefficients
in $\gamma^{25}$. So we conclude that $\Omega^2_5$ is independent from the previous ones. Similar argument
allows to prove that $\Omega^2_6$, $\Omega^2_7$ and $\Omega^2_8$ are also independent from the previous ones.
The form $\Omega^3_4$ is clearly independent from the previous ones. So, at this moment we have 14 independent
curvature forms. Let us consider now the curvature form $\Omega^6_7$. This form could be dependent only of
$\Omega^1_2$ and $\Omega^3_4$. Suppose that $\alpha\,\Omega^1_2 + \beta\, \Omega^3_4 = \frac14, \Omega^6_7$ for
some $\alpha,\beta$. Then, from the coefficients of $\gamma^{58}$ in these curvature forms we get that
$\beta=-\alpha$, but then from the coefficients of $\gamma^{67}$ we conclude that
$\Omega^6_7$ is independent from the previous ones. A similar argument can be applied to $\Omega^6_8$ to get
another independent form.

Therefore there are at least 16 independent curvature forms and this implies that the holonomy is equal to
Spin(7).

\end{document}